\documentclass[journal,10pt]{IEEEtran}
\usepackage{graphicx}
\usepackage{amsfonts}
\usepackage{amssymb}
\usepackage{graphicx}
\usepackage{amsmath}
\usepackage{amsfonts}
\usepackage{amssymb}
\usepackage{epsfig}
\usepackage[named]{algo}
\usepackage{url}

\newcommand{\Section}{\section}
\newcommand{\SubSection}{\subsection}

\newtheorem{theorem}{Theorem}

\newcommand{\bo}[1]{{\bf #1}}

\newcommand{\blambda}{\mbox{\boldmath $\lambda$}}

\title{Multiplicative Noise Removal Using Variable Splitting and Constrained Optimization }
\author{Jos\'{e} M. Bioucas-Dias,~\IEEEmembership{Member,~IEEE}, \hspace{0.5cm}M\'ario A. T. Figueiredo,~\IEEEmembership{Fellow,~IEEE,}
\thanks{J. Bioucas-Dias and M. Figueiredo  are with the {\it Instituto de Telecomunica\c{c}\~{o}es},
   {\it Instituto Superior T\'ecnico,} 1049-001 Lisboa, Portugal; email:
   {\tt mario.figueiredo@lx.it.pt} and {\tt bioucas@lx.it.pt}.}
\thanks{This work was partially supported by {\it Funda\c{c}\~{a}o
para a Ci\^{e}ncia e Tecnologia} (FCT), Portuguese Ministry of Science,
Technology, and Higher Education, under project PTDC/EEA-TEL/104515/2008.}}

\begin{document}

\maketitle

\begin{abstract}
Multiplicative noise (also known as {\it speckle} noise)
models are central to the study of coherent imaging systems,  such as
synthetic aperture radar and sonar, and ultrasound and laser imaging.
These models introduce two additional layers of difficulties
with respect to the standard Gaussian additive noise scenario:
(1) the noise is multiplied by (rather than added to) the original image; (2)
the noise is not Gaussian, with Rayleigh and Gamma being commonly used densities. These
two features of multiplicative noise models preclude the direct application of
most state-of-the-art algorithms, which are designed for solving unconstrained
optimization problems where the objective has two terms: a quadratic data term
(log-likelihood), reflecting the additive and Gaussian nature of the noise,
plus a convex  (possibly nonsmooth)  regularizer ({\em e.g.,} a total
variation or wavelet-based regularizer/prior).  In this paper, we address these
difficulties by: (1)  converting the multiplicative model into
an additive one by taking logarithms, as proposed by some other authors;
(2) using variable splitting to obtain an equivalent constrained problem;
and (3) dealing with  this optimization problem using
the augmented Lagrangian framework. A set of experiments shows that the proposed method,
which we name  MIDAL ({\it multiplicative image denoising by augmented Lagrangian}),
yields state-of-the-art results both in terms of speed and denoising performance.
\end{abstract}



\section{Introduction}
Many special purpose imaging systems use coherent demodulation
of reflected electromagnetic or acoustic waves; well known
examples include ultrasound imaging, synthetic aperture radar (SAR) and
sonar (SAS), and laser imaging. Due to the coherent nature of
these image acquisition
processes, the standard additive Gaussian noise model, so prevalent in
image processing, is inadequate. Instead,
multiplicative noise models, {\it i.e.}, in which the noise field
is multiplied by (not added to) the original image, provide an
accurate description of these coherent imaging systems. In these
models, the noise field is described by a non-Gaussian probability
density function, with Rayleigh and Gamma being common models.

In this introductory section, we begin by briefly recalling how
coherent imaging acquisition processes lead to multiplicative noise
models; for  a more detailed coverage of this topic, the reader is
referred to \cite{conf:Goodman:JOSA:76}, \cite{book:Goodman:GV:07}, \cite{Oliver},
\cite{RudinLionsOsher03}. We then review previous approaches for
dealing with images affected by multiplicative noise and finally
briefly describe the approach proposed in this paper.

\subsection{Coherent Imaging and Multiplicative Noise}
With respect to a given resolution cell of the imaging device, a coherent system acquires
the so-called  {\it in-phase} and {\it quadrature} components, which are the outputs of two
demodulators with respect to, respectively,  $\cos(\omega t)$ and $\sin(\omega t)$, where
$\omega$ is the angular frequency of the carrier signal. Usually, these two components  are collected
into a complex number, with the in-phase and quadrature components corresponding  to the
real and imaginary parts, respectively  \cite{conf:Goodman:JOSA:76}. The complex
observation from a given resolution cell results from the contributions of all the
individual scatterers present in that cell, which interfere in a destructive or
constructive manner, according to their spatial configuration. When this configuration is
random, it yields random fluctuations of the complex observation, a phenomenon which is
usually referred to as {\em speckle}. The statistical properties of speckle have been
widely studied and are the topic of a large body of literature  \cite{conf:Goodman:JOSA:76}, \cite{book:Goodman:GV:07},
\cite{Jakowatz}, \cite{Oliver}.

Assuming a large number of randomly distributed
scatterers and no strong specular reflectors
in each resolution cell  the complex
observation is well modeled by  a zero-mean complex Gaussian circular density ({\em
i.e.}, the real and imaginary parts are independent Gaussian variables with a common
variance). Consequently, the magnitude  of this complex observation follows a Rayleigh
distribution and the square of the magnitude (the so-called intensity) is exponentially
distributed \cite{Jakowatz}, \cite{Oliver}. The term {\it multiplicative noise} is clear
from the following observation: an exponential  random variable can be written as the
product of its mean value, the so-called {\em reflectance} (the parameter of interest to
be estimated) by an exponential variable of unit mean (the noise).

 The scenario just described,  known as {\em fully developed speckle}, leads to observed intensity images with a
characteristic granular appearance,  due to the very low {\em signal to noise ratio} (SNR).
Notice that the SNR, defined as the ratio between the squared intensity mean and the
intensity variance, is equal to one ($0\,$dB); this is a consequence of the equality
between the mean and the standard deviation of an exponential distribution.

\subsection{Improving the SNR: Multi-look Acquisition}
A common approach to improving the SNR in coherent imaging consists in
averaging independent observations of the same resolution cell (pixel). In SAR/SAS systems,
this procedure is called {\em multi-look} ($M$-look, in the case of $M$ looks),
and each independent observation may be  obtained by a different segment of the sensor array.
For fully developed speckle,  the resulting averages are Gamma distributed and
the SNR of an  $M$-look image is improved to $M$.

Another way to obtain an M-look image is to apply a low pass spatial filter (with a
moving average kernel with support size $M$) to a 1-look fully developed speckle image,
making  evident the tradeoff between SNR and spatial resolution. This type of $M$-look
image can be understood as an estimate of the underlying reflectance, under the
assumption that this reflectance is constant in the set of cells (pixels) included in the
averaging process. A great deal of research has been devoted to developing space variant
filters which average  large numbers of pixels in homogeneous regions, yet avoid
smoothing across reflectance discontinuities, in order to preserve  details/edges
\cite{art:FrostStilesPAMI82}. Many other speckle reduction techniques have been proposed;
see \cite{Oliver} for a comprehensive review of the literature up to 1998.

\subsection{Estimation of Reflectance: Variational Approaches}
\label{sec:estimation_var}
What is usually referred to as multiplicative noise removal is of course nothing
but the estimation of the reflectance of the underlying scene. This is an inverse
problem calling for regularization, which usually consists in assuming
 that the underlying reflectance image
is piecewise smooth. In image denoising under multiplicative noise, this assumption has
been formalized, in a Bayesian estimation framework, using  Markov random field
priors \cite{Bioucas-Dias98}, \cite{Oliver}. More recently, variational approaches
using {\em total variation} (TV) regularization were proposed
\cite{AubertAujol08}, \cite{Huang09}, \cite{RudinLionsOsher03}, \cite{ShiOsher07},
\cite{Steidl09}.  In a way, these approaches extend the
spatial filtering methods referred to in the previous subsection; instead of  explicitly
piecewise flat estimates, these approaches yield piecewise smooth estimates adapted to
the structure of the underlying reflectance.

Both the variational and the Bayesian {\it maximum a posteriori} (MAP) formulations to image
denoising (under multiplicative, Gaussian, or other noise models) lead to optimization
problems with two terms: a data fidelity term (log-likelihood) and a regularizer (log-prior).
Whereas in Gaussian noise models, the data fidelity term is quadratic, thus quite benign from an optimization
point of view, the same is no longer true under multiplicative observations.  In
\cite{AubertAujol08},  the  data fidelity term is the negative log-likelihood resulting
directly from  the $M$-look multiplicative model, which, being non-convex,  raises
difficulties from an optimization point of view. Another class of approaches, which
is the one adopted in this paper, also uses the $M$-look multiplicative model, but
yields  convex data fidelity terms by formulating the problem with respect to the
logarithm of the  reflectance; see \cite{DurandFadiliN08}, \cite{Huang09},
\cite{RudinLionsOsher03}, \cite{ShiOsher07}, \cite{Steidl09}, and references therein.
A detailed analysis of several data fidelity terms for the  multiplicative noise
model can be found in \cite{Steidl09}.

The combination of TV regularization with the log-likelihood resulting from the
multiplicative  observation model  leads to an objective function, with a non-quadratic
term (the log-likelihood) plus a non-smooth term (the TV regularizer), to which some
research work has been recently devoted \cite{AubertAujol08}, \cite{Huang09},
\cite{RudinLionsOsher03}, \cite{ShiOsher07}, \cite{Steidl09}.  Even when the
log-likelihood is convex (as in \cite{Huang09}, \cite{RudinLionsOsher03},
\cite{ShiOsher07}, \cite{Steidl09}), it does not have a Lipschitz-continuous gradient,
which is a necessary condition for the applicability (with guaranteed convergence) of
algorithms of the forward-backward splitting (FBS) class \cite{fista}, \cite{twist},
\cite{CombettesSIAM},  \cite{sparsa}.  Methods based on the Douglas-Rachford  splitting
(DRS), which do not require the objective function to have a Lipschitz-continuous
gradient, have been recently proposed \cite{combettes2007}, \cite{Steidl09}.

\subsection{Proposed Approach}

\label{sec:sub_approach}

In this paper we address the (unconstrained) convex optimization problem which results
from the $M$-look multiplicative model formulated with respect to the logarithm of the
reflectance. As shown in \cite{Steidl09}, this is the most adequate formulation to
address reflectance estimation under multiplicative noise with TV regularization. We
propose an optimization algorithm with the following building blocks:
\begin{itemize}
  \item the original unconstrained optimization problem is first transformed into an equivalent
        constrained problem, via a variable splitting procedure;
  \item this constrained problem is then addressed using an augmented
        Lagrangian method.
\end{itemize}

This paper is an elaboration of  our previous work \cite{BioucasFigueiredoICIP09}, where
we have  addressed multiplicative noise removal also via a variable splitting procedure.
In that paper, the  constrained optimization problem was dealt with by exploiting the recent
{\it split-Bregman} approach \cite{GoldsteinOsher}, but using a splitting strategy which is quite
different from the one in \cite{GoldsteinOsher}.

It happens that the Bregman iterative methods recently proposed to handle imaging inverse
problems are equivalent to {\it augmented Lagrangian} (AL) methods \cite{Nocedal}, as shown in
\cite{YinOsherGoldfarbDarbon}. We prefer the AL perspective, rather than the Bregman iterative
view, as it is a standard and more elementary optimization tool (covered in most textbooks on
optimization). In particular, we solve the constrained problem resulting from the variable
splitting using an algorithm (of the AL family) known as {\it alternating direction method
of multipliers} (ADMM)  \cite{EcksteinBertsekas}, \cite{Gabay}, \cite{Glowinski}.

Other authors have recently addressed the variational restoration of speckled images
 \cite{AubertAujol08}, \cite{combettes2007}, \cite{Huang09},
\cite{RudinLionsOsher03}, \cite{ShiOsher07}, \cite{Steidl09}. The commonalities and
differences between our approach and the approaches followed by other authors
will be discussed after the detailed description of our method, since this discussion
requires notation and concepts which will be introduced in the next section.\\

The paper is organized as follows. Section \ref{sec:formulation} formulates the problem,
including the detailed description of the multiplicative noise model and the TV
regularization we adopt to estimate the reflectance image.
Section \ref{sec:basic_tools} reviews the  variable splitting and the
augmented Lagrangian optimization methods, which are the basic tools with which our
approach is built. Section \ref{sec:approach}  introduces the proposed algorithm by
direct application of the basic tools introduced in Section \ref{sec:basic_tools}.
Section~\ref{sec:related} discusses related work.
Section \ref{sec:experiments} reports the results of a series of experiments aiming at
comparing the proposed algorithm with previous state-of-the-art competitors.
Finally, Section \ref{sec:conclusion} concludes the paper.

\section{Problem Formulation}
\label{sec:formulation} Let $\bo{y}\in \mathbb{R}_{+}^{n}$ denote an $n$-pixel observed
image, assumed to be a sample of a random image $\bo{ Y}$, the mean of which is the
underlying (unknown) reflectance image $\bo{x}\in \mathbb{R}_{+}^{n}$, that is,
$\mathbb{E}[\bo{Y}] = \bo{x}$. Adopting a conditionally independent multiplicative noise
model, we have
\begin{equation}
Y_i = x_i  \, N_i ,\;\; \mbox{for\  $\; i=1,...,n$,}\label{eq:multiply}
\end{equation}
where ${\bf N}\in \mathbb{R}_{+}^{n}$ is an image of independent and identically distributed (iid)
random variables with unit mean, $\mathbb{E}(N_i)=1$, following a common probability density function $p_N$.
In the case of $M$-look fully developed speckle noise, $p_N$ is a Gamma density
\begin{equation}
   \label{eq:Gamma}
   p_N(n) = \frac{M^M}{\Gamma(M)}\; n^{M-1}e^{-nM},
\end{equation}
thus with expected value  $\mathbb{E}[N]=1$ and  variance
\[
 \sigma_N^2\equiv \mathbb{E}\left[(N-\mathbb{E}[N])^2\right]=\frac{1}{M}.
\]
Accordingly, we define the signal-to-noise ratio (SNR) associated to a random
variable $Y_i$, for $i=1,\dots,n$, as
\begin{equation}\label{eq:snr}
    \mbox{SNR}\equiv \frac{\mathbb{E}[Y_i]^2}{\sigma^2_{Y_i}} = M.
\end{equation}

An additive noise model is obtained by taking logarithms of (\ref{eq:multiply})
\cite{DurandFadiliN08}, \cite{Huang09},
\cite{RudinLionsOsher03}, \cite{ShiOsher07}, \cite{Steidl09}.
For an arbitrary pixel of the image (thus dropping the pixel subscript for simplicity),
the observation model becomes
\begin{equation}
      \label{eq:log_observation}
       \underbrace{\log Y}_G = \underbrace{\log x}_z + \underbrace{\log N}_W.
\end{equation}
The density of the random variable $W = \log N$ is
\begin{eqnarray}
      p_W(w) & =  & p_N(e^w)\,e^w \nonumber\\
             & =  & \frac{M^M}{\Gamma(M)}\; e^{Mw}e^{-e^w M},
\end{eqnarray}
thus
\begin{equation}
    p_{G|z}(g|z) = p_W(g-z).
\end{equation}
Invoking  the conditional independence assumption, we are finally lead to
\begin{eqnarray}
\log p_{\bo{G}|\bo{z}}(\bo{g}|\bo{z}) & = & \sum_{s=1}^n \log p_{W}(g_s - z_s)\\
& = & C - M \sum_{s=1}^n (z_s + e^{\,g_s-z_s}) ,
\end{eqnarray}
where $C$ is a constant (irrelevant for estimation purposes).

Using the MAP criterion (which is equivalent to a regularization method),
the original image is inferred by solving an unconstrained minimization problem with the form
\begin{equation}
       \widehat{\bo{z}}   \in \arg\min_{\bo{z}} L(\bo{z}),  \label{eq:L_uncons}
\end{equation}
where $L(\bo{z})$ is the objective function given by
\begin{eqnarray}
L(\bo{z}) & = & -\log p_{\bo{G}|\bo{Z}}(\bo{g}|\bo{z}) +\lambda\, \phi(\bo{z}).\\
& = & M \sum_{s=1}^n  \left(z_s + e^{g_s-z_s}\right) + \lambda\, \phi(\bo{z}) + A;\label{eq:neg_like}
\end{eqnarray}
in (\ref{eq:neg_like}), $A$ is an irrelevant additive constant,  $\phi$ is the
negative of the log-prior (the regularizer), and $\lambda$ is the so-called
regularization parameter.

In this  work, we adopt the standard isotropic discrete TV regularizer \cite{Chambolle04}, that is,
\begin{equation}
 \phi(\bo{z})  \equiv \sum_{s=1}^n \sqrt{(\Delta^h_s\bo{z})^2+(\Delta^v_s\bo{z})^2},\label{eq:theTV}
\end{equation}
where  $\Delta^h_s\bo{z}$ and $\Delta^v_s\bo{z}$ denote the horizontal and vertical first
order differences at pixel $s\in\{1,\dots,n\}$, respectively.

Each term $\left(z_s + e^{g_s-z_s} \right)$ of (\ref{eq:neg_like}),  corresponding to the
negative log-likelihood, is strictly convex and coercive, thus so is their sum. Since the
TV regularizer is also convex (though not strictly so), the objective function $L$
possesses a unique minimizer \cite{CombettesSIAM}, which is a fundamental property, in
terms of optimization. In contrast,  the formulation of the problem in terms of the
original variables $\bo{x}$ (rather than their logarithm) leads to a non-convex
optimization problem \cite{AubertAujol08}, \cite{RudinLionsOsher03}. As seen in
\cite{AubertAujol08}, the uniqueness of the minimizer of that non-convex objective
is not guaranteed in general.

In this paper, we will address problem (\ref{eq:L_uncons}) using   variable splitting and
augmented Lagrangian optimization. In the next section, we briefly review these tools,
before presenting our approach in detail.

\Section{Basic Tools} \label{sec:basic_tools}

\SubSection{Variable Splitting}
Consider an unconstrained optimization problem in which
the objective function is the sum of two functions, one
of which is written as  a composition:
\begin{equation}
\min_{{\bf u}\in \mathbb{R}^n} f_1({\bf u}) + f_2\left(g({\bf u})\right).\label{unconstrained_basic}
\end{equation}

Variable splitting is a very simple procedure that consists in
creating a new variable, say ${\bf v}$,
to serve as the argument of $f_2$, under the
constraint that $g({\bf u}) = {\bf v}$. The idea is to
consider the constrained problem
\begin{equation}\begin{array}{cl}
 {\displaystyle \min_{{\bf u},{\bf v}\in \mathbb{R}^n}} & f_1({\bf u}) + f_2({\bf v})\\
 \mbox{s.t.} &g({\bf u}) = {\bf v}, \end{array}\label{constrained_basic}
\end{equation}
which is clearly equivalent to unconstrained problem (\ref{unconstrained_basic}).
Notice that in the feasible set $\{({\bf u},{\bf v}): g({\bf u}) = {\bf v}\}$, the objective
function in (\ref{constrained_basic}) coincides with that in (\ref{unconstrained_basic}).

The rationale behind variable splitting methods is that it may
be easier to solve the constrained problem (\ref{constrained_basic})
than it is to solve its unconstrained counterpart (\ref{unconstrained_basic}).
This idea has been recently used in several image processing applications
\cite{BioucasFigueiredo2008}, \cite{GoldsteinOsher}, \cite{Huang09}, \cite{WangYangYinZhang}.

A variable splitting method was used in \cite{WangYangYinZhang} to obtain a fast
algorithm for TV-based image restoration. Variable splitting was also used in
\cite{BioucasFigueiredo2008} to handle problems involving compound regularizers; {\it
i.e.}, where instead of a single regularizer $\lambda \phi$ as in
(\ref{eq:neg_like}), one
has a linear combination of two (or more) regularizers $\lambda_1 \phi_1 + \lambda_2
\phi_2$. In \cite{BioucasFigueiredo2008}, \cite{Huang09}, and \cite{WangYangYinZhang},
the constrained problem (\ref{constrained_basic}) is attacked by a quadratic penalty
approach, {\it i.e.}, by solving
\begin{equation}
 \min_{{\bf u},{\bf v}\in \mathbb{R}^n}  f_1({\bf u}) + f_2({\bf v}) +
 \frac{\mu}{2}\, \|g({\bf u}) -{\bf v}\|_2^2,
\label{quadratic_penalty}
\end{equation}
by alternating minimization with respect to ${\bf u}$ and ${\bf v}$, while
slowly taking $\mu$ to very large values (a {\it continuation} process),
to force the solution of (\ref{quadratic_penalty}) to approach that
of (\ref{constrained_basic}), which in turn is equivalent to (\ref{unconstrained_basic}).
The rationale behind these methods is that each step of this
alternating minimization may be much easier than the original
unconstrained problem (\ref{unconstrained_basic}). The drawback is that
as $\mu$ becomes very large, the intermediate minimization problems
become increasingly ill-conditioned, thus causing numerical problems
(see \cite{Nocedal}, Chapter 17).

A similar variable splitting approach underlies the recently
proposed split-Bregman methods \cite{GoldsteinOsher}; however,
instead of using a quadratic penalty technique, those methods
attack the constrained problem directly using a Bregman iterative
algorithm \cite{YinOsherGoldfarbDarbon}. It has been shown that,
when $g$ is a linear function, {\it i.e.}, the constraints in
(\ref{constrained_basic}) are linear ({\it e.g.}, $g({\bf u}) = {\bf Gu}$),
the Bregman iterative
algorithm is equivalent to the augmented Lagrangian method
\cite{YinOsherGoldfarbDarbon}, which is briefly reviewed in the
following subsection.

\SubSection{Augmented Lagrangian}
In this brief presentation of the augmented Lagrangian method,
we closely follow \cite{Nocedal}, which the reader should
consult for more details.
Consider an equality constrained optimization problem (which
includes (\ref{constrained_basic}) as a particular instance,
if $g$ is linear)
\begin{equation}\begin{array}{cl}
 {\displaystyle \min_{{\bf z}\in \mathbb{R}^d}} & E({\bf z})\\
 \mbox{s.t.} & {\bf A z - b = }\mbox{\boldmath $0$}, \end{array}\label{constrained_linear}
\end{equation}
where ${\bf b} \in \mathbb{R}^p$ and ${\bf A}\in \mathbb{R}^{p\times d}$,
{\it i.e.}, there are $p$ linear equality constraints.
The so-called augmented Lagrangian function for this problem is
defined as
\begin{equation}
{\cal L}_A ({\bf z},\blambda,\mu) = E({\bf z}) + \blambda^T ({\bf Az-b}) +
\frac{\mu}{2}\,  \|{\bf Az-b}\|_2^2,\label{augmented_L}
\end{equation}
where $\blambda \in \mathbb{R}^p$ is a vector of Lagrange multipliers
and $\mu \geq 0$ is the penalty parameter.

The so-called {\it augmented Lagrangian  method}
(ALM), also known as the {\it method of multipliers}
\cite{Hestenes}, \cite{Powell},  consists in minimizing ${\cal L}_A ({\bf
z},\blambda,\mu)$ with respect to ${\bf z}$, keeping $\blambda$ fixed, and then updating
$\blambda$.

\vspace{0.3cm}
\begin{algorithm}{ALM}{\label{alg:salsa1}}
Set $k=0$, choose $\mu > 0$, and  $\blambda_0$.\\
\qrepeat\\
      ${\bf z}_{k+1} \in \arg\min_{{\bf z}} {\cal L}_A ({\bf z},\blambda_k,\mu)$\\
     $\blambda_{k+1} \leftarrow \blambda_{k} + \mu ({\bf Az}_{k+1} - {\bf b})$\\
     $k \leftarrow k + 1$
\quntil stopping criterion is satisfied.
\end{algorithm}
\vspace{0.3cm}

Although it is possible (even recommended) to update the value of $\mu$ in each iteration
\cite{Bazaraa}, \cite{Nocedal}, we will not consider that  option in this paper.
Importantly, unlike in the quadratic penalty method, it is not necessary to take $\mu$
to infinity to guarantee that the ALM converges to the solution of the
constrained problem (\ref{constrained_linear}).

Notice that (after a straightforward complete-the-squares procedure) the terms added to
$E({\bf z})$ in the definition of the augmented Lagrangian  ${\cal L}_A ({\bf
z},\blambda,\mu)$ in (\ref{augmented_L}) can be written as a single quadratic term,
leading to the following alternative form for the ALM algorithm:

\vspace{0.3cm}
\begin{algorithm}{ALM (version II)}{\label{alg:salsa2}}
Set $k=0$, choose $\mu > 0$, ${\bf z}_0$, and  ${\bf d}_0$.\\
\qrepeat\\
     ${\bf z}_{k+1} \in \arg\min_{{\bf z}} E({\bf z}) + \frac{\mu}{2}\|{\bf Az-d}_k\|_2^2$\\
     ${\bf d}_{k+1} \leftarrow {\bf d}_{k} - ({\bf Az}_{k+1} - {\bf b})$\\
     $k \leftarrow k+1$
\quntil stopping criterion is satisfied.
\end{algorithm}
\vspace{0.3cm}

This form of the ALM algorithm makes clear its equivalence with the Bregman iterative
method (see \cite{YinOsherGoldfarbDarbon}).

It has been shown that, with adequate initializations, the ALM algorithm generates the
same sequence as a {\it proximal point algorithm} (PPA) applied to the Lagrange dual of
problem (\ref{constrained_linear}); for further details, see \cite{Iusem}, \cite{Setzer},
and references therein. Moreover, the sequence $\{{\bf d}_{k}\}$
converges to a solution of this dual problem and that all cluster points of the sequence
$\{{\bf z}_{k}\}$ are solutions of the (primal) problem (\ref{constrained_linear})
\cite{Iusem}.

\SubSection{Augmented Lagrangian for Variable Splitting} We now show how ALM can be used
to address problem (\ref{constrained_basic}), in the particular case where $g({\bf u}) \equiv {\bf u}$
 ({\it i.e.}, $g$ is the identity function).
 This problem can be written in the form (\ref{constrained_linear}) using
the following definitions:
\begin{equation}
{\bf z} \equiv \left[\!\begin{array}{c}{\bf u}\\{\bf v}\end{array}\!\right],\hspace{0.5cm}
{\bf b} \equiv {\bf 0}, \hspace{0.5cm}
{\bf A} \equiv  [\,{\bf I} \;\; -{\bf I}\,],
\end{equation}
and
\begin{equation}
E({\bf z}) \equiv f_1({\bf u}) + f_2({\bf v}).
\end{equation}
With these definitions in place, steps 3 and 4 of the ALM (version II) can be written as
follows
\begin{eqnarray}
\left[\!\!\!\begin{array}{c}{\bf u}_{k+1} \\ {\bf v}_{k+1} \end{array}\!\!\!\right] \!\!  & \!\!  \in \!\! & \!\!  \arg\min_{{\bf u},{\bf v}}
  \left\{\! f_{1}({\bf u}) \!  + \!
  f_{2}({\bf v})  \! + \! \frac{\mu}{2} \|{\bf u} \! -\! {\bf v} \! - \!{\bf d}_k\|_2^2
 \right\} \hspace{0.4cm}\label{mixed}\\
{\bf d}_{k+1}\!\!  & \!\!  \leftarrow \!\! & \!\!  {\bf d}_{k} - ( {\bf u}_{k+1} - {\bf v}_{k+1}).
\end{eqnarray}

The minimization (\ref{mixed}) is not trivial since, in general, it involves
non-separable quadratic and possibly non-smooth terms. A natural solution is to use a
{\it non-linear block-Gauss-Seidel} (NLBGS) technique, in which (\ref{mixed}) is solved by
alternating minimization with respect to ${\bf u}$ and ${\bf v}$.
Of course this raises several questions: for a given ${\bf d}_k$,
how much computational effort should be spent in approximating the solution of
(\ref{mixed})? Does this NLBGS procedure converge?
Taking just one step of this NLBGS scheme in each iteration of ALM
leads to an algorithm known as the {\it alternating direction method of multipliers}
(ADMM)  \cite{EcksteinBertsekas}, \cite{Gabay}, \cite{Glowinski} (see also \cite{Esser},
\cite{Setzer}, \cite{Steidl09}):

\vspace{0.3cm}
\begin{algorithm}{ADMM}{
\label{alg:salsa3}}
Set $k=0$, choose $\mu > 0$, ${\bf v}_0$, and  ${\bf d}_0$.\\
\qrepeat\\
   $  {\bf u}_{k+1}  \in  \arg\min_{{\bf u}} f_{1}({\bf u})
 + \frac{\mu}{2} \|{\bf u} - {\bf v}_k - {\bf d}_k\|_2^2$\\
  $  {\bf v}_{k+1}  \in  \arg\min_{{\bf v}} f_{2}({\bf v})
 + \frac{\mu}{2} \|{\bf u}_{k+1} - {\bf v} - {\bf d}_k\|_2^2$\\
     ${\bf d}_{k+1} \leftarrow {\bf d}_{k} - ({\bf u}_{k+1} - {\bf v}_{k+1} - {\bf b})$\\
     $k \leftarrow k+1$
\quntil stopping criterion is satisfied.
\end{algorithm}
\vspace{0.3cm}

For later reference, we now recall the theorem by Eckstein and Bertsekas
\cite{EcksteinBertsekas} in which convergence of (a generalized version of) ADMM is
shown. This theorem applies to problems of the form (\ref{unconstrained_basic}) with
$g({\bf u}) \equiv {\bf G}{\bf u}$, {\it i.e.},
\begin{equation}
\min_{{\bf u}\in \mathbb{R}^n} f_1({\bf u}) + f_2\left({\bf G\,
u}\right),\label{unconstrained_basic_linear}
\end{equation}
of which (\ref{constrained_basic})  is the constrained optimization
reformulation.

\vspace{0.3cm}
\begin{theorem}[Eckstein-Bertsekas, \cite{EcksteinBertsekas}]
\label{th:Eckstein}{\sl Consider problem  (\ref{unconstrained_basic_linear}), where
 ${\bf G}\in\mathbb{R}^{d\times n}$ has full column rank and
 $f_1\,:\,\mathbb{R}^n\rightarrow (-\infty,\infty]$ and $f_2\,:\,\mathbb{R}^d\rightarrow
 (-\infty,\infty]$ are closed, proper, convex functions. Consider arbitrary
 $\mu>0$ and ${\bf v}_0, {\bf d}_0\in \mathbb{R}^d$.
 Let $\{\eta_k \geq 0, \; k=0,1,...\}$ and $\{\nu_k \geq 0,
\; k=0,1,...\}$ be two sequences such that
\[
\sum_{k=0}^\infty \eta_k < \infty \;\;\;\mbox{and} \;\;\; \sum_{k=0}^\infty \nu_k <
\infty.
\]
Consider three sequences $\{{\bf u}_k \in \mathbb{R}^{n}, \; k=0,1,...\}$, $\{{\bf v}_k
\in \mathbb{R}^{d}, \; k=0,1,...\}$, and $\{{\bf d}_k \in \mathbb{R}^{d}, \; k=0,1,...\}$
that satisfy
\begin{eqnarray}
 \eta_k & \geq & \left\| {\bf u}_{k+1} - \arg\min_{{\bf u}} f_{1}({\bf u})
 + \frac{\mu}{2} \|{\bf G}{\bf u} \! - \!{\bf v}_k \! -\! {\bf d}_k\|_2^2 \right\| \nonumber\\
 \nu_k & \geq & \left\| {\bf v}_{k+1}  - \arg\min_{{\bf v}} f_{2}({\bf v})
 + \frac{\mu}{2} \|{\bf G}{\bf u}_{k+1} \! - \! {\bf v} \! - \! {\bf d}_k\|_2^2 \right\| \nonumber\\
 {\bf d}_{k+1} & = & {\bf d}_{k} - ({\bf Gu}_{k+1} - {\bf v}_{k+1}).\nonumber
 \end{eqnarray}
Then, if (\ref{unconstrained_basic_linear}) has a solution, say ${\bf u}^*$,
the sequence $\{{\bf u}_k\}$ to ${\bf u}^*$. If (\ref{unconstrained_basic_linear}) does not have a
solution, then at least one of the sequences $\{{\bf v}_k \}$ or $\{{\bf d}_k\}$
diverges.}
\end{theorem}
\vspace{0.3cm}

Notice that the ADMM algorithm defined above generates sequences $\{{\bf u}_k\}$, $\{{\bf
v}_k\}$, and $\{{\bf d}_k \}$ which satisfy the conditions in Theorem \ref{th:Eckstein}
in a strict sense ({\it i.e.}, with $\eta_k = \nu_k = 0$). One of the important
consequences of this theorem is that it shows that it is not necessary to exactly solve
the minimizations in lines 3 and 4 of ADMM; as long as the sequence of errors are
absolutely summable, convergence is not compromised.

The proof of Theorem \ref{th:Eckstein} is based on the equivalence between ADMM and the
so-called Douglas-Rachford splitting method (DRSM) applied to the dual of  problem
(\ref{unconstrained_basic_linear}). The DRSM was recently used for image recovery
problems in \cite{combettes2007}. For recent and comprehensive reviews of ALM, ADMM,
DRSM, and their relationship with Bregman and split-Bregman methods, see \cite{Esser},
\cite{Setzer}.

\section{Proposed Approach}
\label{sec:approach} To address the optimization problem  (\ref{eq:L_uncons}) using the
tools reviewed  in the previous section, we begin by rewriting  it as
 \begin{eqnarray}
      \label{eq:L_cons}
      (\widehat{\bo{z}},\widehat{\bo{u}}) & = &
\arg\min_{\bo{z},\bo{u}} L(\bo{z},\bo{u})\\
      \label{eq:z_u_c}
      \mbox{s.t.} & &   \bo{z} = \bo{u},
\end{eqnarray}
 with
\begin{eqnarray}
      \label{eq:L_cons_obj}
     L(\bo{z},\bo{u})  = M \sum_{s=1}^n \left( z_s + e^{g_s-z_s}\right) + \lambda\,\phi(\bo{u}).
\end{eqnarray}
Notice how the original variable (image) $\bo{z}$ was split into a pair of variables
$(\bo{z,u})$,  which are decoupled in the objective function (\ref{eq:L_cons_obj}).

The approach followed in \cite{Steidl09} also considers a variable splitting, aiming at
the application of the ADMM method. However, the splitting therein adopted is different
from ours and, as shown below, leads to a more complicated algorithm with an additional
ADMM inner loop.

Applying the ADMM method to the constrained problem defined by (\ref{eq:L_cons})--(\ref{eq:L_cons_obj})  leads to the proposed algorithm, which we call {\it multiplicative
image denoising by augmented Lagrangian} (MIDAL). Obviously,  the estimate of the image
${\bf x}$ is computed as  $\widehat{\bf x} = e^{{\bf z}_k}$, component-wise.

\vspace{0.3cm}
\begin{algorithm}{MIDAL}{
\label{alg:midal}}
Choose $\mu > 0$, $\lambda \geq 0$, ${\bf u}_0$, and  ${\bf d}_0$. Set $k \leftarrow  0$.\\
\qrepeat\\
    $\bo{z}' \leftarrow \bo{u}_k +\bo{d}_k$\\
   $\bo{z}_{k+1}  \leftarrow {\displaystyle \arg\min_{\bo{z}}  \sum_{s=1}^n \left( z_s +
e^{\, g_s-z_s}\right)+\frac{\mu}{2M}\| \bo{z}-\bo{z}'\|^2}$\\
     $\bo{u}' \leftarrow \bo{z}_k-\bo{d}_k$\\
  $\bo{u}_{k+1} \leftarrow  {\displaystyle \arg\min_{\bo{u}} \frac{1}{2}\|
\bo{u}-\bo{u}'\|^2+\frac{\lambda}{\mu}\,\phi(\bo{u})}$.\\
     ${\bf d}_{k+1} \leftarrow {\bf d}_{k} - ({\bf z}_{k+1} - {\bf u}_{k+1})$\\
     $k \leftarrow k+1$
\quntil a stopping criterion is satisfied.
\end{algorithm}
\vspace{0.3cm}

The minimization with respect to $\bo{z}$ (line 4) is in fact a set of $n$ decoupled scalar convex
minimizations. Each of these minimizations has closed form solution in terms of the
Lambert W function \cite{Corless}. However, as in \cite{BioucasFigueiredoICIP09}, we
apply the Newton method, as it  yields a faster  (and very accurate) solution by running just
a few iterations.

The minimization  with respect to $\bo{u}$ (line 6) corresponds to solving a $\ell_2-$TV denoising
problem with observed image ${\bf u}'$ and regularization parameter $\lambda/\mu$ or,
equivalently, to  computing the so-called Moreau proximity operator (see \cite{CombettesSIAM})
of $(\lambda/\mu)\phi$, denoted $\mbox{prox}_{(\lambda/\mu)\phi}$ at ${\bf u}'$; {\em i.e.}, for
$\gamma \geq 0$,
\begin{equation}
    \label{eq:prox_TV}
     \mbox{prox}_{\gamma\phi}({\bf u}')\equiv \arg\min_{\bo{u}} \frac{1}{2}\|
\bo{u}-\bo{u}'\|^2+\gamma\phi(\bo{u}).
\end{equation}
We use Chambolle's algorithm \cite{Chambolle04} to compute $\mbox{prox}_{\gamma\phi}$,
although faster algorithms could be applied \cite{WangYangYinZhang}.  As stated in
Theorem \ref{th:Eckstein}, this computation does not have to be solved exactly as long as
the Euclidian norm of the  errors are  summable  along the ADMM iterations (and thus along
the MIDAL iterations).

Still invoking Theorem \ref{th:Eckstein}, and assuming that the  sequences of
optimization errors with respect to ${\bf z}$ (line 4 of MIDAL pseudo-code) and  ${\bf
u}$ (line 6 of MIDAL pseudo-code) are absolutely summable, then MIDAL convergence is
guaranteed, because $f_1$ and $f_2$ are closed proper convex functions, and ${\bf
G} = \bf I$ has full column rank.

\section{Comments on Related Work}
\label{sec:related}
We will now make a few remarks on related work. Notice how the variable
splitting approach followed by the ADMM method allowed converting
a  difficult problem involving a non-quadratic term and a TV regularizer into
a sequence of two simpler problems: a decoupled minimization problem and a
TV denoising problem with a quadratic data term. In contrast, the variable
splitting adopted in \cite{Steidl09}  leads to an intermediate  optimization
that is neither separable nor quadratic, which is dealt with
by an inner DRS iterative technique.

TV-based image restoration under multiplicative noise was recently addressed in
\cite{ShiOsher07}. The authors apply an inverse scale space flow, which
converges to the solution of the constrained problem of minimizing $\mbox{TV}({\bf z})$
under an equality constraint on the log-likelihood; this requires a carefully
chosen stopping criterion, because the solution of this constrained problem
is not a good estimate. Moreover, as evidenced in the experiments reported
in \cite{DurandFadiliN08}, the method proposed in \cite{ShiOsher07}
has a performance far from the state-of-the-art.

In \cite{Huang09}, a  variable splitting  is also used to obtain
an objective function with the form
\begin{equation}
E({\bf z,u}) = L({\bf z,u}) + \alpha \|{\bf z - u}\|_2^2;
\end{equation}
this is the so-called splitting-and-penalty method. Notice that the
minimizers of $E({\bf z,u})$ converge to those of (\ref{eq:L_cons})-(\ref{eq:z_u_c})
only when $\alpha$ approaches infinity. However, since
$E({\bf z,u})$ becomes severely ill-conditioned when $\alpha$ is very
large, causing numerical difficulties, it is only practical to minimize
$E({\bf z,u})$ with moderate values of $\alpha$; consequently, the solutions obtained
are not minima of the regularized negative log-likelihood (\ref{eq:neg_like}).
Nevertheless, the method exhibits good performance, although not as
good as the method herein proposed, as shown in the experiments reported below.

As mentioned in Subsection~\ref{sec:estimation_var}, in the approach followed in \cite{AubertAujol08},  the  objective function is con-convex. In addition to
a lack of guarantee of uniqueness of the minimizer, this feature raises
difficulties from an optimization point of view. Namely, as confirmed
experimentally in \cite{Huang09}, the obtained estimate depends
critically on the initialization.

Finally, we should mention that the algorithmic approach herein pursued
can also be interpreted from a Douglas-Rachford splitting perspective
\cite{combettes2007}, \cite{Esser}, \cite{Setzer}. In \cite{combettes2007},
that approach was applied to several image restoration problems with
non-Gaussian noise, including a multiplicative noise case, but not with the
Gamma distribution herein considered.

\section{Experiments}
\label{sec:experiments}
In this section, we report experimental results comparing the
performance of the proposed approach with those of the recent state-of-the-art methods
introduced in \cite{DurandFadiliN08} and \cite{Huang09}. We chose to focus on those
two references for two reasons: (a) they report quantitative results and
the corresponding implementations are available; (b) experimental results reported
in those papers show that the methods therein proposed outperform other recent techniques,
namely the above mentioned \cite{AubertAujol08} and \cite{ShiOsher07}, as well as the
(non-iterative) block-Stein thresholding of \cite{Chesnau_Fadili_Starck}.

The proposed algorithm is implemented in MATLAB 7.5 and all the tests were carried out on
a PC with a 3.0GHz Intel Core2Extreme CPU and 4Gb of RAM. All the
experiments use synthetic data, in the sense that the observed image is generated
according to (\ref{eq:multiply})--(\ref{eq:Gamma}), where ${\bf x}$ is the original image.
In Table \ref{tab:expsetup} we list the details of the 16 experimental setups
considered: the original images, their sizes, and the maximum and
minimum pixel values ($x_{\mbox{\footnotesize max}}$ and $x_{\mbox{\footnotesize min}}$);
the M values (which coincides with the SNR (\ref{eq:snr})); the adopted value of
$\lambda$ (the regularization parameter in (\ref{eq:neg_like})) for our algorithm.

\begin{figure}
\centering
\includegraphics[width=0.48\columnwidth,height=0.48\columnwidth]{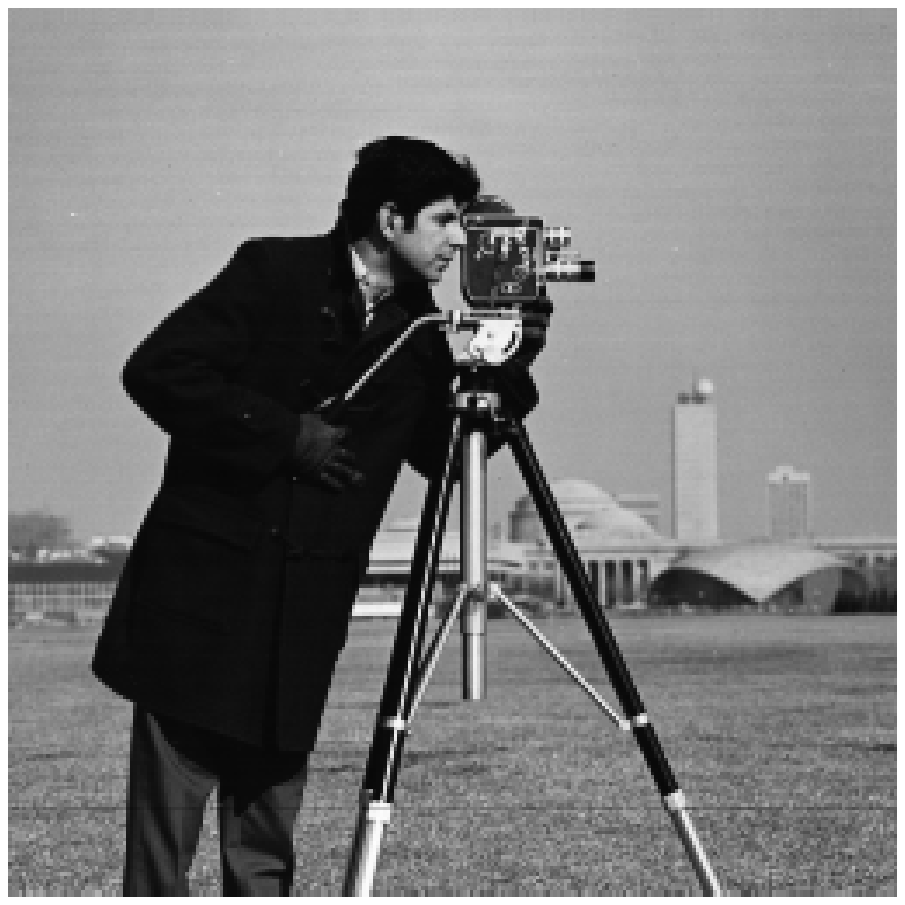}\hspace{.019\columnwidth}\includegraphics[width=0.48\columnwidth,height=0.48\columnwidth]{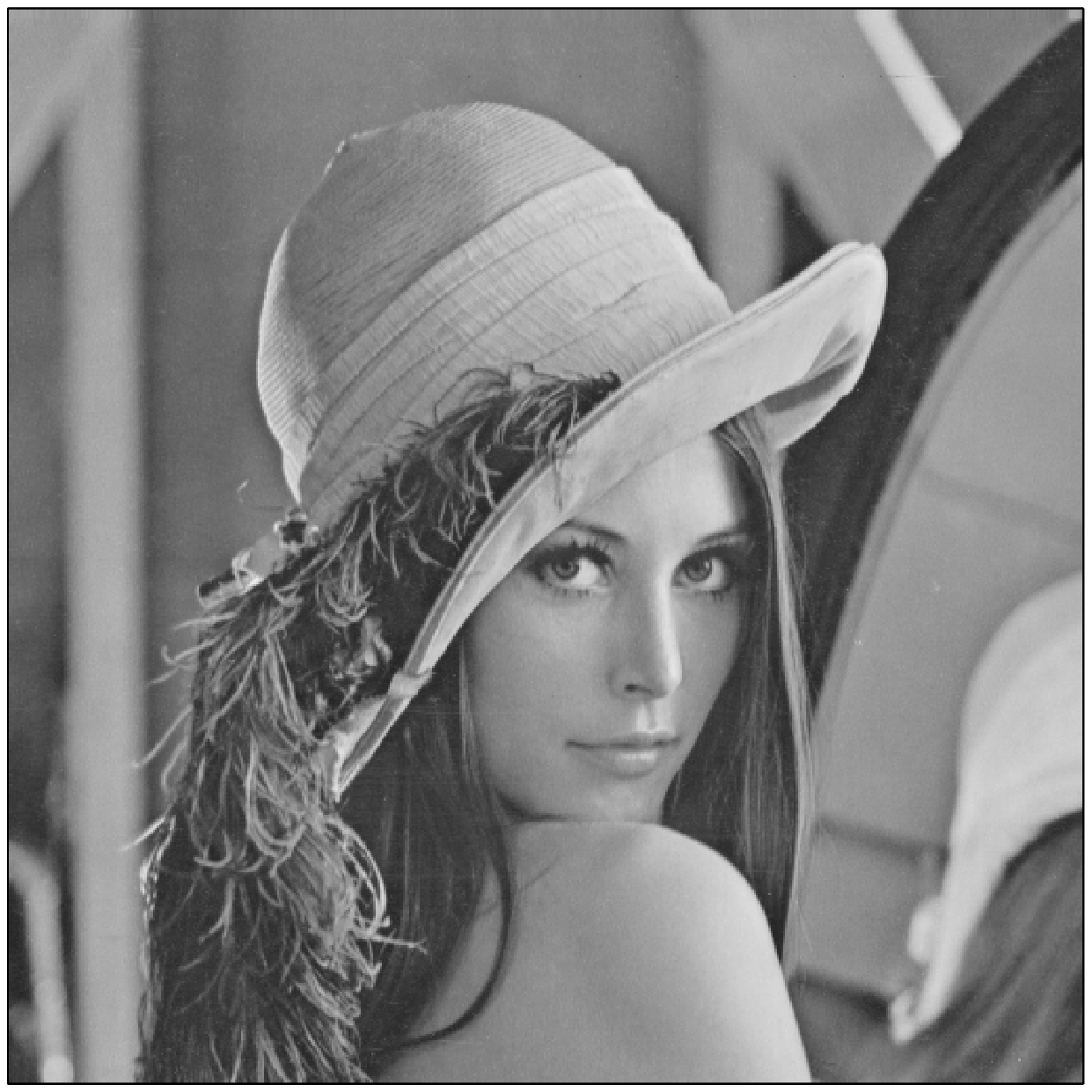}

\vspace{0.1cm}
\includegraphics[width=0.48\columnwidth,height=0.48\columnwidth]{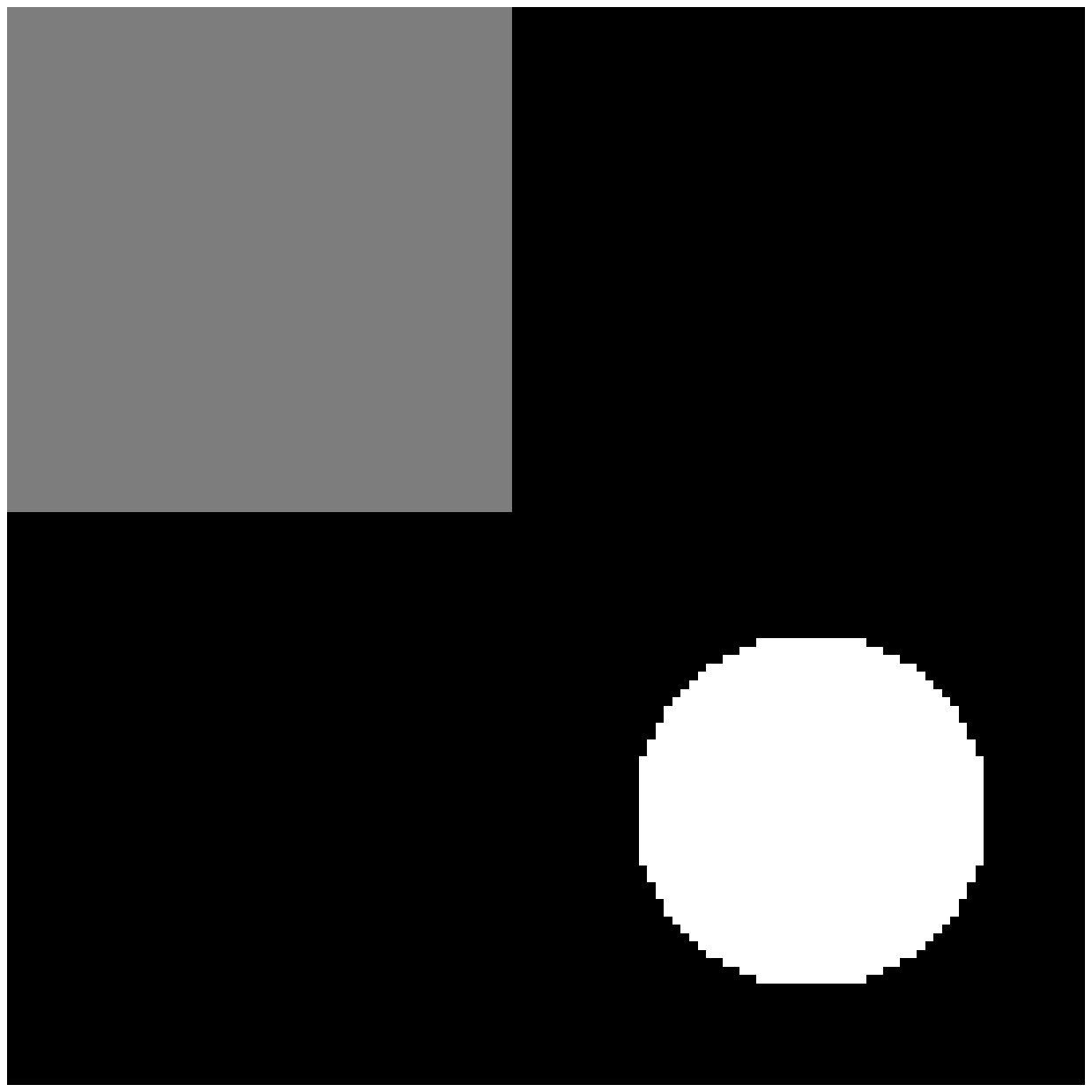}\hspace{.019\columnwidth}\includegraphics[width=0.48\columnwidth,height=0.48\columnwidth]{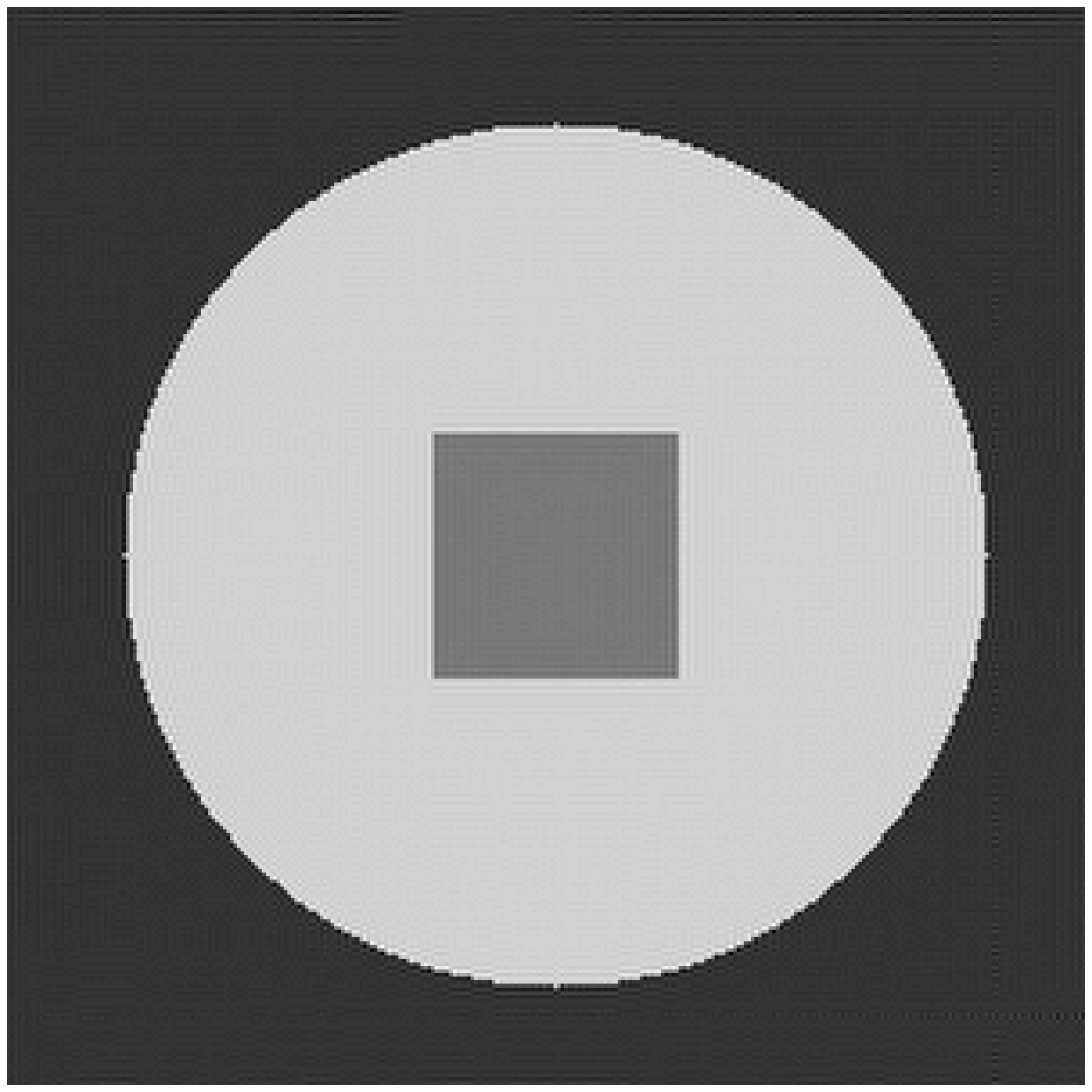}

\vspace{0.1cm}
\includegraphics[width=0.48\columnwidth,height=0.48\columnwidth]{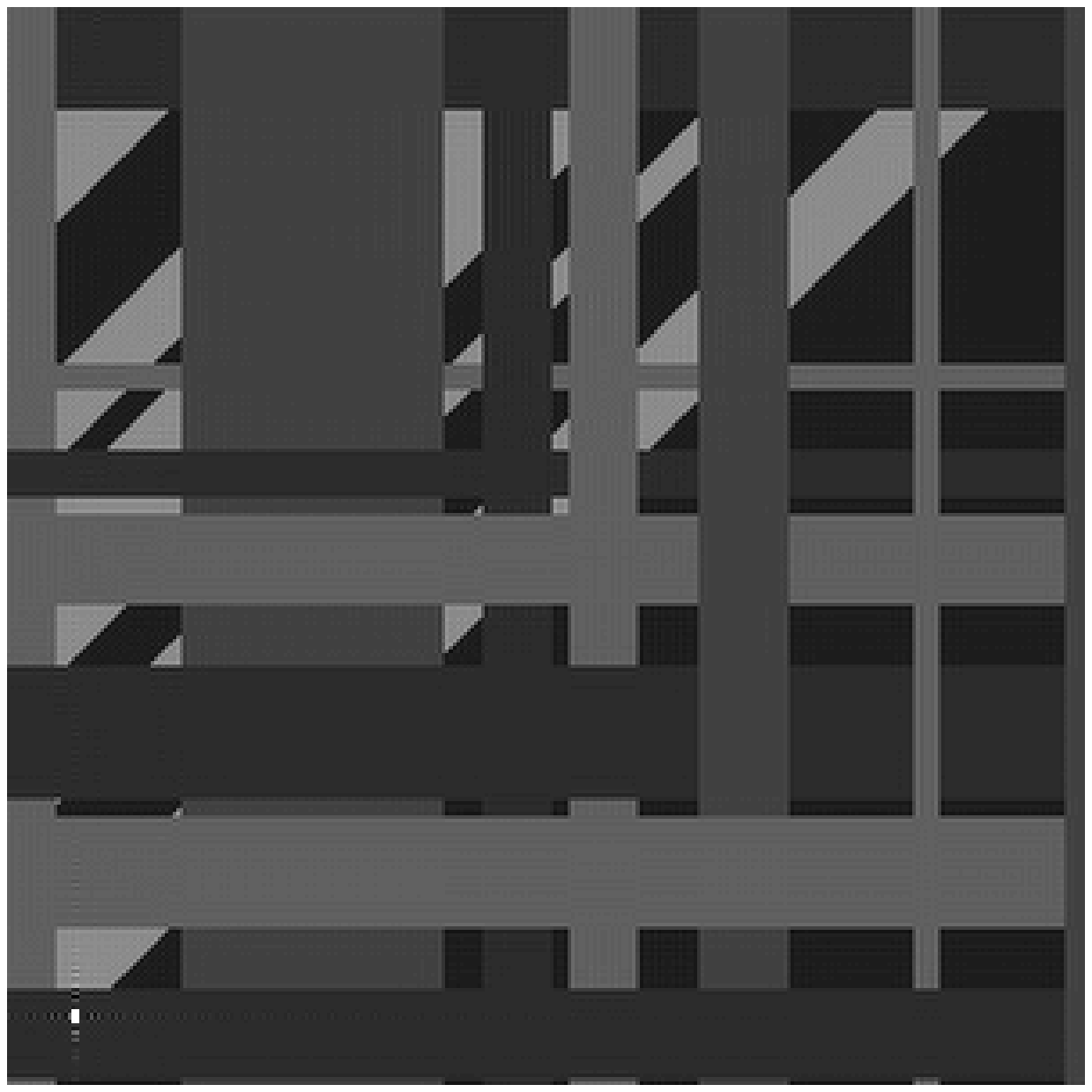}\hspace{.019\columnwidth}\includegraphics[width=0.48\columnwidth,height=0.48\columnwidth]{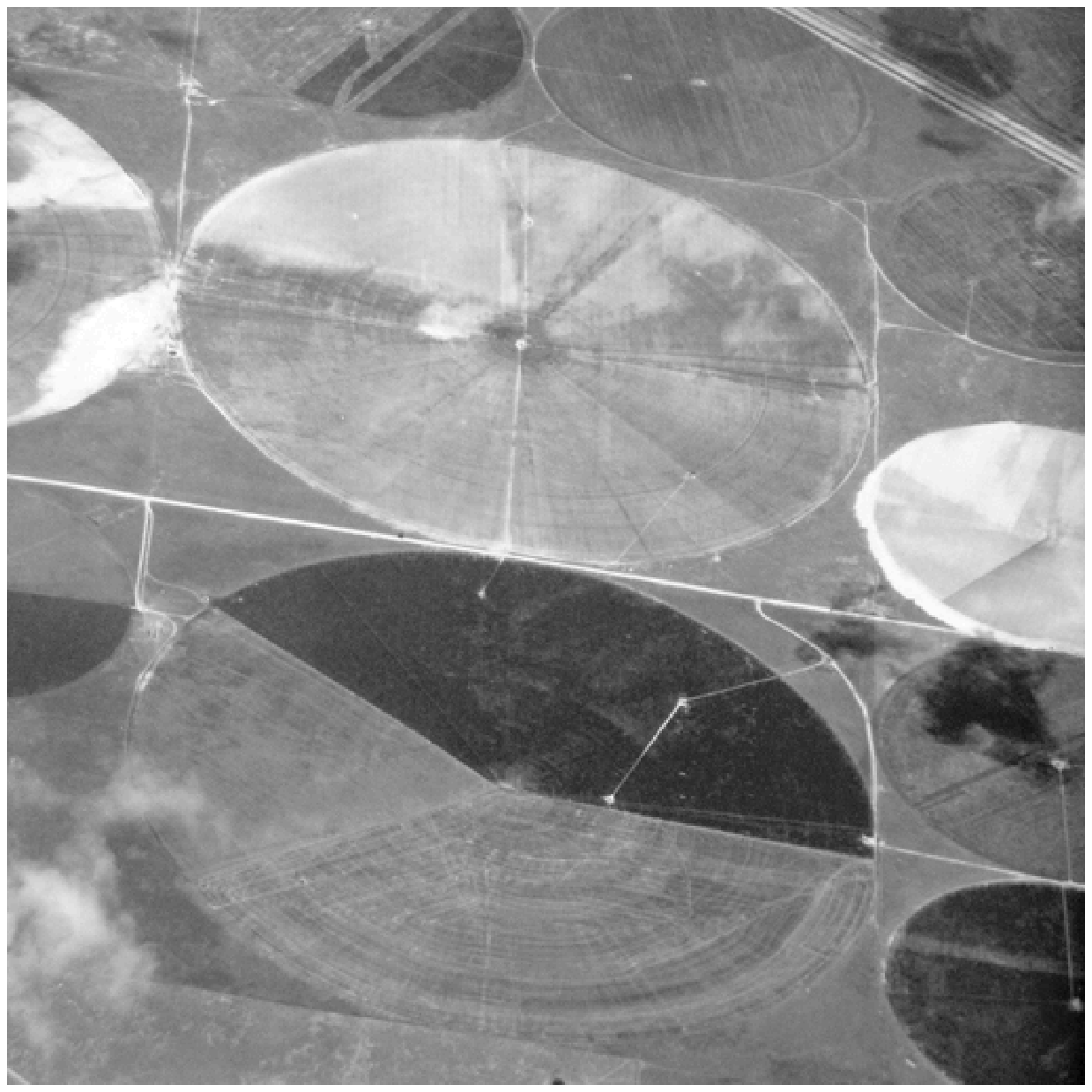}

\vspace{0.1cm}
\centerline{\includegraphics[width=0.48\columnwidth,height=0.48\columnwidth]{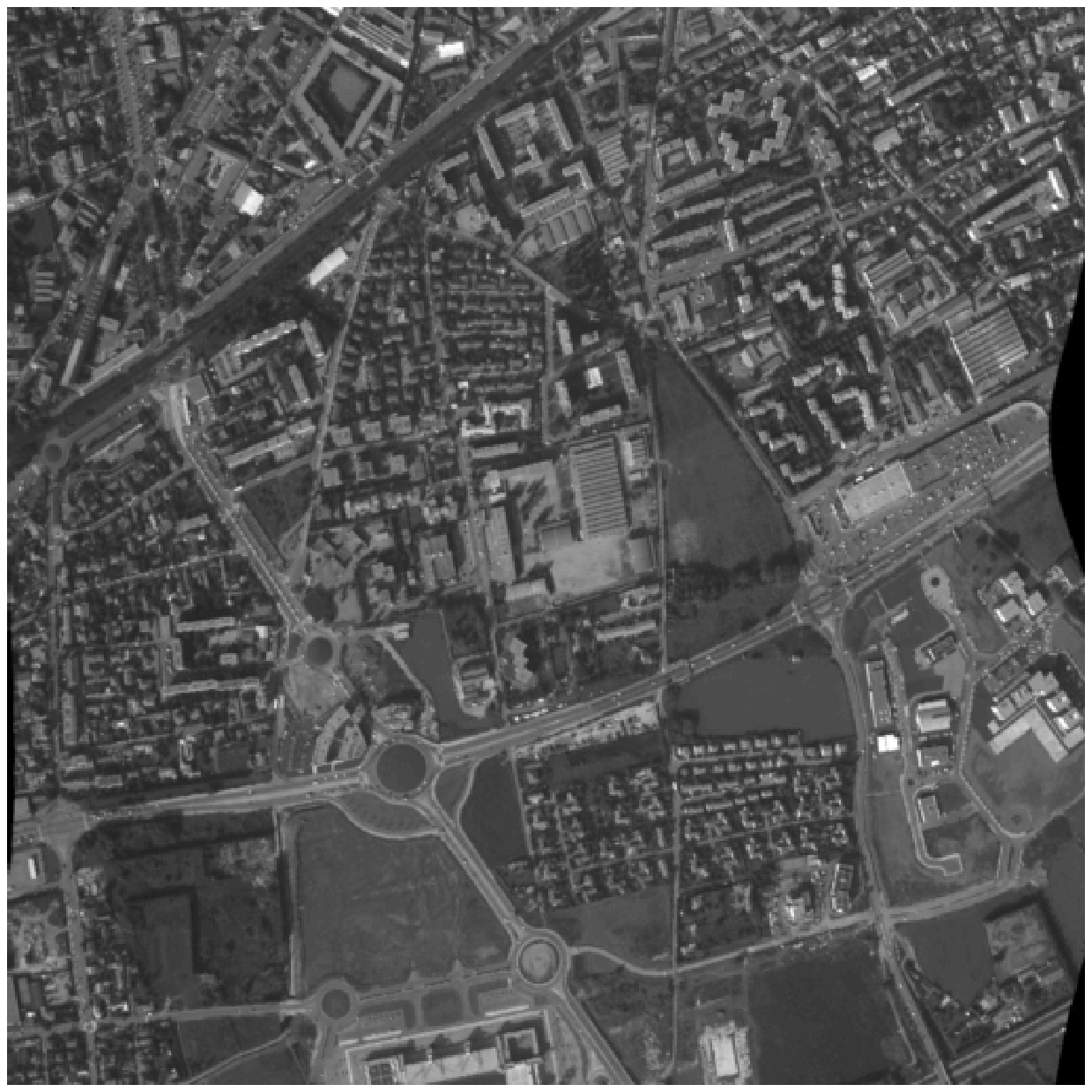}}

\caption{The 7 original images used in the 16 experiments: from top to bottom and
left to right: Cameraman, Lena, Sim1, Sim2, Sim3, Fields, and N\^{\i}mes.} \label{fig:originals}
\end{figure}

Experiments $1-7$ reproduce the experimental setup
used in \cite{AubertAujol08} and \cite{Huang09}; experiments 8--16 follow those
reported in \cite{DurandFadiliN08}. The 8 original images used are shown in
Figure \ref{fig:originals}. The values of $x_{\mbox{\footnotesize max}}$ and
$x_{\mbox{\footnotesize min}}$ are as in \cite{DurandFadiliN08} and \cite{Huang09}, for
comparison purposes. Notice the very low  SNR values ($M$ values) for most  observed
images, a usual scenario in applications involving multiplicative noise.

The focus of this paper is mainly the speed of the algorithms to solve the optimization
problem (\ref{eq:L_uncons}), thus the automatic choice of the regularization parameter
$\lambda$ is out of scope. Therefore,  as in \cite{AubertAujol08},
\cite{DurandFadiliN08},  and \cite{Huang09},  we select $\lambda$ by searching for the
value leading to the lowest mean squared error with respect to the true image.

Assuming  that conditions of Theorem \ref{th:Eckstein} are met,  MIDAL is guaranteed to
converge for any value of the penalty parameter $\mu>0$.  This parameter has, however, a
strong impact in the convergence speed. We have verified experimentally that setting
$\mu = \lambda$ yields good
results. For these reason, we have used this setting in all the experiments.

MIDAL  is initialized with the observed noisy image. The quality of the estimates is
assessed using  the  relative error (as in \cite{Huang09})
\[
\mbox{Err} \equiv \frac{\|\widehat{\bf x} - {\bf x}\|_2}{\|{\bf x}\|_2},
\]
and the mean absolute-deviation error (MAE) (as in  \cite{DurandFadiliN08})
\[
   \mbox{MAE} \equiv   \frac{\|\widehat{\bf x} - {\bf x}\|_1}{n},
\]
where $\widehat{\bf x}\equiv \exp(\widehat{\bf z})$ and  $\|\cdot \|_2$ and $\|\cdot \|_1$ stand for the $\ell_2$ and $\ell_1$ norms, respectively. As in
\cite{Huang09}, we use the stopping criterion
\[
     \frac{\|{\bf x}_{k+1} - {\bf x}_{k}\|_2}{\|{\bf x}_{k}\|_2}\leq 10^{-m},
\]
with $m=4$ in experiments 1 to 7, as in \cite{Huang09}, and $m=2$ in experiments 8 to 16.

\begin{table}
\centering \caption{Experimental setup}\label{tab:expsetup} \vspace{0.1cm}
\begin{tabular}{ | c |  l | c | c | c | c | c |}
\hline & & & & & & \\
   Exp.  &  Image   & Size      & $M$ & $x_{\mbox{\footnotesize min}} $ & $x_{\mbox{\footnotesize max}}$ & $\lambda$\\
    & & & & & & \\ \hline  \hline
    1 & Cameraman & $256\times 256$ & 3     & 0.03 & 0.9 & 4\\
    2 & Cameraman & $256\times 256$ & 13    & 0.03 & 0.9 & 6.5\\
    3 & Lena      & $256\times 256$      & 5     & 0.03 & 0.9 & 4.5\\
    4 & Lena      & $256\times 256$      & 33    & 0.03 & 0.9 & 9.8\\
    5 & Sim1      & $128\times 128$      & 2     & 0.05 & 0.46 & 5.5 \\
    6 & Sim2      & $300\times 300$      & 1     & 0.13 & 0.39 & 3.5 \\
    7 & Sim3      & $300\times 300$      & 2     & 0.03 & 0.9 & 3\\
    8 & Fields     & $512\times 512$      & 1    & 101 & 727 & 3.5\\
    9 & Fields      & $512\times 512$      & 4    & 101 & 727 & 4.5\\
    10 & Fields  & $512\times 512$      & 10    & 101 & 727 & 6.7\\
    11 & N\^imes     & $512\times 512$      & 1    & $10^{-4}$ & 255 & 2\\
    12 & N\^imes      & $512\times 512$      & 4    & $10^{-4}$ & 255 & 2.7\\
    13 & N\^imes  & $512\times 512$      & 10    & $10^{-4}$ & 255 & 4\\
    14 & Cameraman     & $256\times 256$      & 1    & 7 & 253 & 2.7\\
    15 & Cameraman      & $256\times 256$      & 4    & 7 & 253 & 4.5\\
    16 & Cameraman  & $256\times 256$      & 10    & 7 & 253 & 6.1\\
\hline
\end{tabular}
\end{table}

\subsection{Computing the TV Proximity Operator}
MIDAL requires, in each iteration,  the computation of the TV proximity
operator, $\mbox{prox}_{\gamma\phi}$, given by (\ref{eq:prox_TV}), for which
we use Chambolle's fixed point iterative algorithm \cite{Chambolle04}.
Aiming at faster convergence of Chambolle's algorithm, and
consequently of MIDAL, we initialize each run with the dual variables (see \cite{Chambolle04}
for details on the dual variables) computed in the previous run.
The underlying rationale is that, as MIDAL proceeds, the images to which the
proximity operator is applied get closer; thus, by initializing the computation
of the next proximity operator with the internal variables of the previous iteration,
the burn-in period is largely avoided.

\begin{figure}
\centerline{\includegraphics[width=0.8\columnwidth]{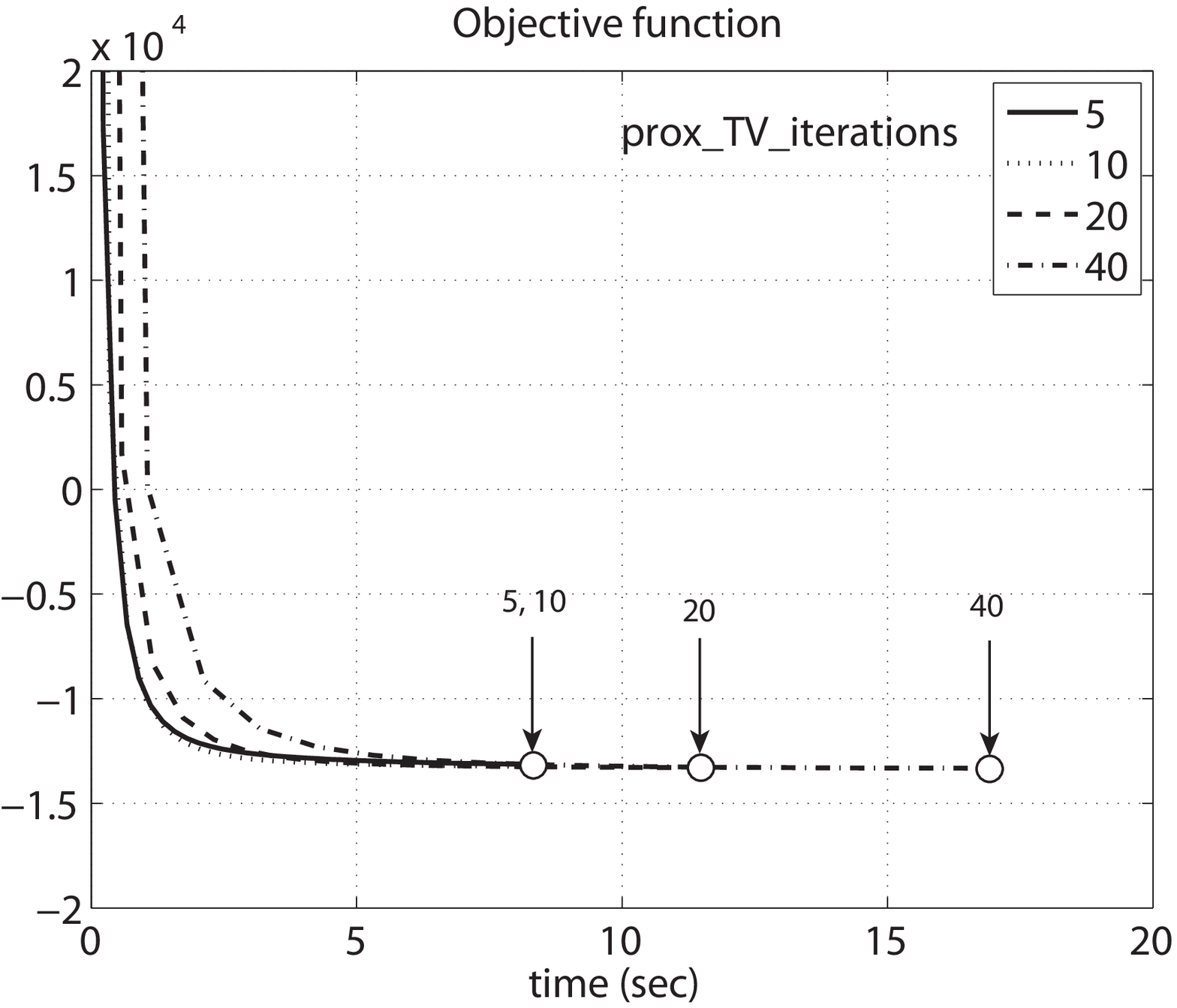}}

\centerline{\includegraphics[width=0.8\columnwidth]{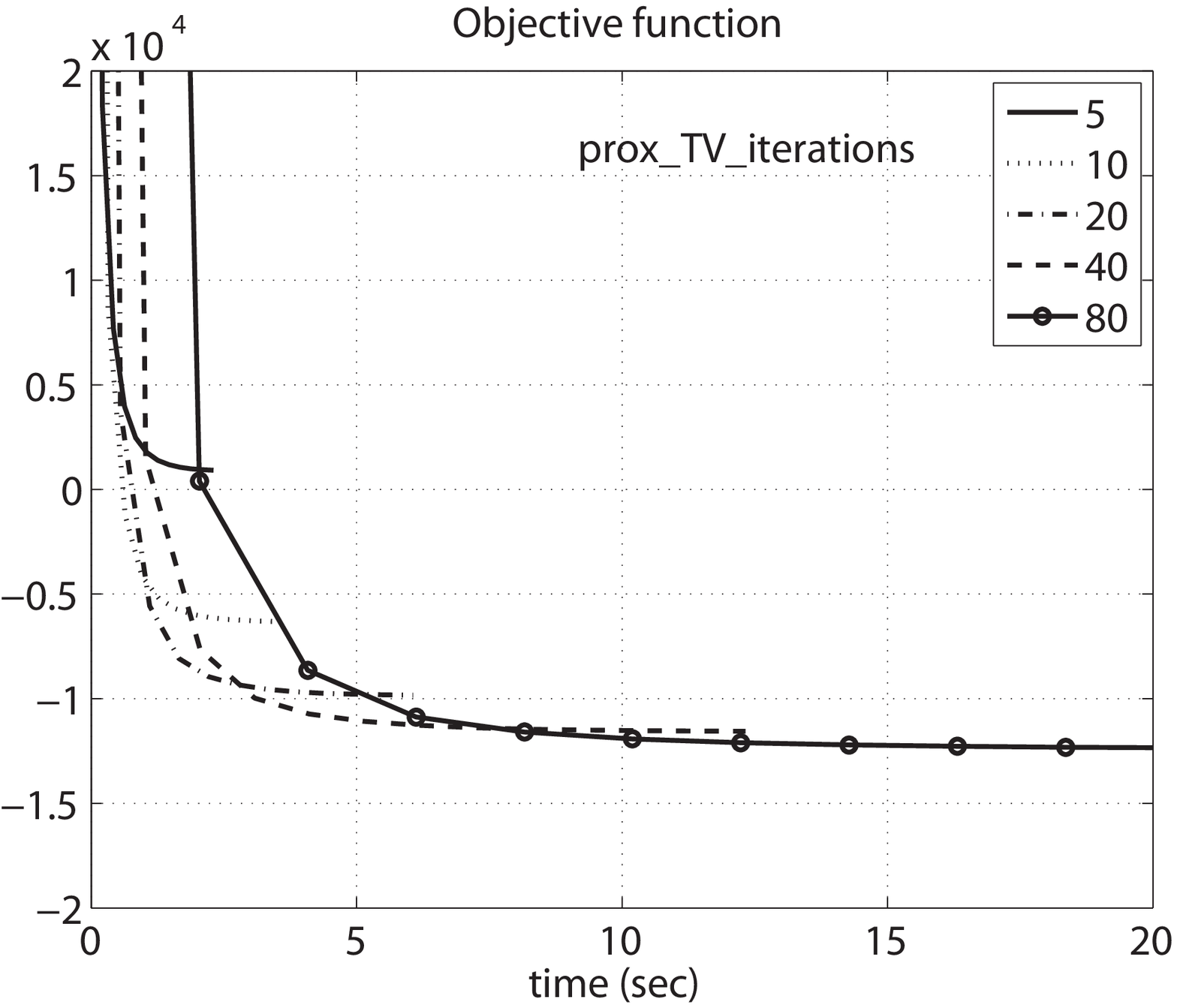}}

\caption{Evolution of the objective function (\ref{eq:neg_like}) for the setting of
Experiment 1, using Chambolle's fixed point  iterative algorithm \cite{Chambolle04}
to compute the TV proximity operator with different initializations. Top plot:
initialization with the dual variables computed in the previous iteration.
Bottom plot: initialization  with the dual variables set to zero.} \label{fig:dual_vars}
\end{figure}

Another perspective to look at this procedure, already referred to, is given by Theorem
\ref{th:Eckstein}, which states that there is no need to exactly solve the minimizations
in each iteration, but just to ensure the minimization errors along the iterations are
absolutely summable. The fulfilment of this condition is easier to  achieve with the proposed
initialization  than with a fixed initialization. Fig. \ref{fig:dual_vars} illustrates
this aspect. For the setting of Experiment 1,  it shows the evolution of the objective
function (\ref{eq:neg_like}) when the dual variables are initialized  with the ones
computed in the previous iteration (left hand) and when the dual variables are
initialized to zero (right hand).  All the curves in the top plot reach essentially
the same value. Notice that  MIDAL  takes, approximately, the same time for a number of fixed
point iterations between 5 and 20 to compute the TV proximity operator. For a number of
iterations higher than 20, MIDAL time increases because, in each iteration, it runs  more
fixed point  iterations than necessary. In the plots in the right hand side, we see
that the minimum of the objective function is never reached, although it can be
approximated for large values of the fixed point iterations. Based on these observations,
we set the number of  fixed point iterations to 20 in all experiments of this section.

\subsection{Results}

Table~\ref{tab:results} reports the results of the 16  experiments. The times for
our algorithm and that of \cite{Huang09} are relative to the computer mentioned above.
The numbers of iterations are also given, but just as side information since the
computational complexity per iteration of each algorithm is different. The
initialization of the algorithm of \cite{Huang09} is either the observed image
or the mean of the observed image; since the final values of Err are essentially
the same for  both initializations, we report the best of the two times.

For experiments 8 to 16, we did not have access to the code of \cite{DurandFadiliN08},
so we report the MAE and Err values presented in that paper. According to the authors,
their algorithm was run for a fixed number of iterations, thus the computation time
depends only on the image size. The time values shown in Table~\ref{tab:results}
were provided by the authors and were obtained on a MacBook Pro with a 2.53GHz Intel CoreDuo
processor and 4Gb of RAM.

In  experiments 1 to 7, our method always achieves lower estimation errors than
the method of \cite{Huang09}. Notice that the gain of MIDAL is  larger for images
with lower SNR, corresponding to the more difficult problems. Moreover, our
algorithm is faster than that of \cite{Huang09} in all the experiments, by a
factor larger than 3.

In all the experiments 8 to 16, our algorithm achieves lower MAE than the
method of \cite{DurandFadiliN08}. Concerning the relative error Err, our
algorithm outperforms theirs in 5 out of 9 cases, there is a tie in two
cases, and is outperformed (albeit by a very small margin) in two cases.

\begin{table*}
\centering \caption{Experimental results. Iter, Err, and MAE denote, respectively,  the
number of iterations, the relative error, and the mean absolute-deviation
error. The times are reported in seconds. The time for \cite{DurandFadiliN08} is referred to a
different machine (see text).}\label{tab:results} \vspace{0.3cm}
\begin{tabular}{c | l c l l | l l l| l l l }
\hline\hline
& \multicolumn{4}{c|}{MIDAL} & \multicolumn{3}{c|}{\cite{Huang09}} &  \multicolumn{3}{c}{\cite{DurandFadiliN08}}\\
Exp. & Err  & MAE   & Iter  & Time & Err  & Iter & Time &  Err & MAE & Time$^*$\\ \hline
 1      &  \textbf{0.130} & 0.035 & 21 & \textbf{10} & 0.151 & 113 & 32 & -- & -- & --\\
\hline
 2      &  \textbf{0.090} & 0.025 & 16 & \textbf{8} & 0.098 & 115 & 36 &  -- & -- & --\\
\hline
 3      &  \textbf{0.111} & 0.036 & 23 & \textbf{10} & 0.118 & 133 & 38 &  -- & -- & --\\
\hline
 4      &  \textbf{0.069} & 0.023 & 17 & \textbf{9} & 0.071 & 182 & 52 &  --  & -- & --\\
\hline
 5      &  \textbf{0.128} & 0.009 & 16 & \textbf{1.4} & 0.143 & 165 & 8 &  --  & -- & --\\
\hline
 6      &  \textbf{0.069} & 0.012 & 32 & \textbf{21} & 0.083 & 166 & 70 &  --  & -- & --\\
\hline
 7      &  \textbf{0.137} & 0.023 & 17 & \textbf{11} & 0.174 & 165 & 69 & -- & -- & -- \\
\hline
 8     &  \textbf{0.089} & \textbf{29.06 }    & 61 &  121  & -- & -- & -- &  0.096 & 32.67 & 245 \\
\hline
 9     &  \textbf{0.066} & \textbf{20.86 }   & 34 &  68  & -- & -- & -- &  \textbf{0.066} & 22.0 & 245\\
\hline
 10     &   0.056 & \textbf{17.94 }& 29 &  58  & -- & -- & -- &  \textbf{0.055}  & 18.24 & 245 \\
\hline
 11     &  \textbf{0.301} & \textbf{12.38 }    & 25 &  14  & -- & -- & -- &  0.314 & 13.27 & 245 \\
\hline
 12     &  \textbf{0.217} & \textbf{8.78 }    & 25 &  50  & -- & -- & -- &  \textbf{0.217}  & 8.98 & 245\\
\hline
 13     &  \textbf{0.170} & \textbf{6.87  }   & 23 &  46 & -- & -- & -- &  0.174 & 7.11 & 245 \\
\hline
 14     &  \textbf{0.167} & \textbf{12.74}     & 33 &  15  & -- & -- & -- &  0.192 & 16.78 & 54 \\
\hline
 15     &  \textbf{0.124} & \textbf{9.43}    & 19 &  9  & -- & -- & -- &  0.131 & 10.67 & 54 \\
\hline
 16     &  0.097 & \textbf{7.42}  & 56 &  26  & -- & -- & -- &  \textbf{0.091} & 7.44 & 54 \\
\hline
\end{tabular}
\end{table*}

\begin{figure}
\centering
\includegraphics[width=0.48\columnwidth,height=0.48\columnwidth]{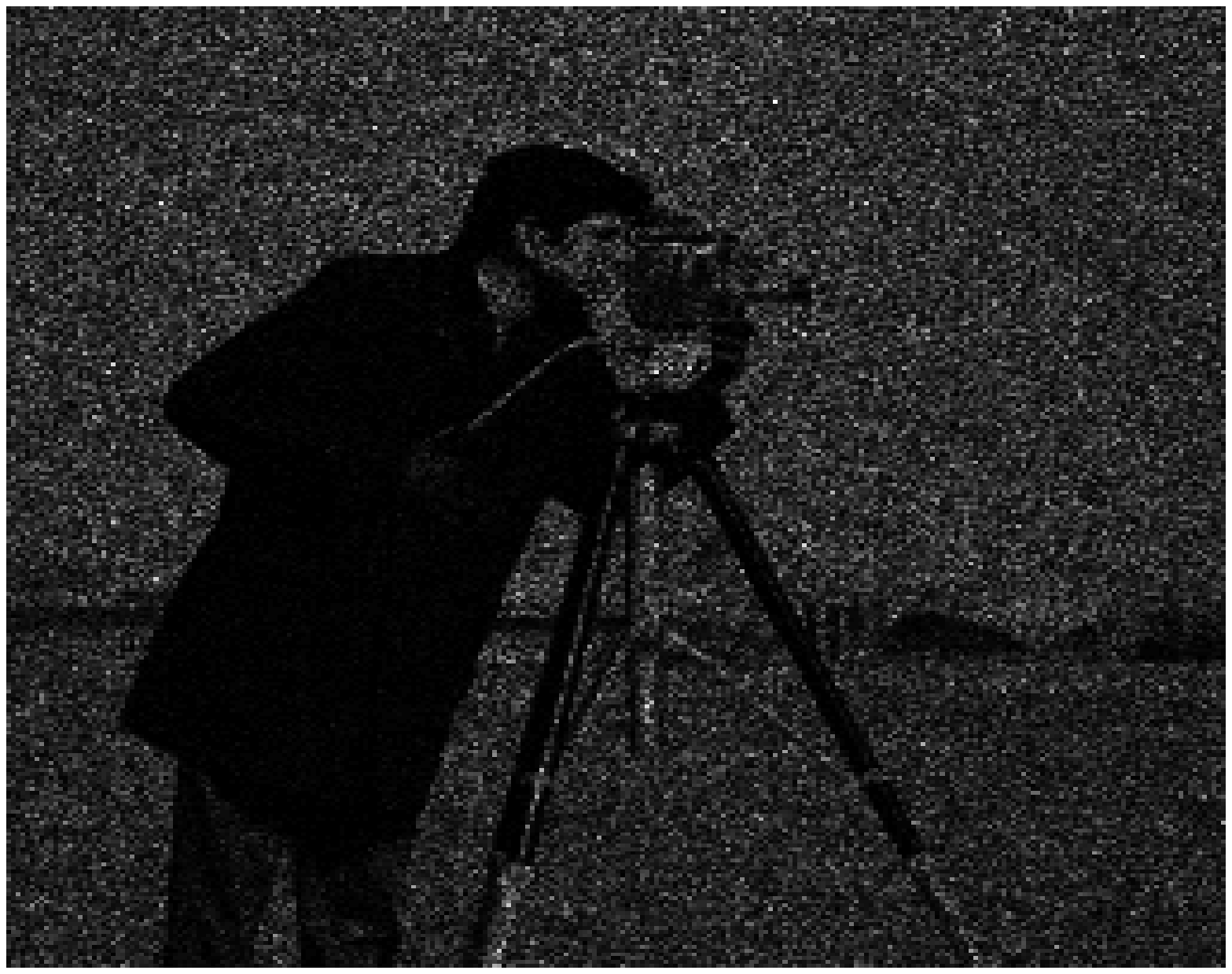}\hspace{.019\columnwidth}\includegraphics[width=0.48\columnwidth,height=0.48\columnwidth]{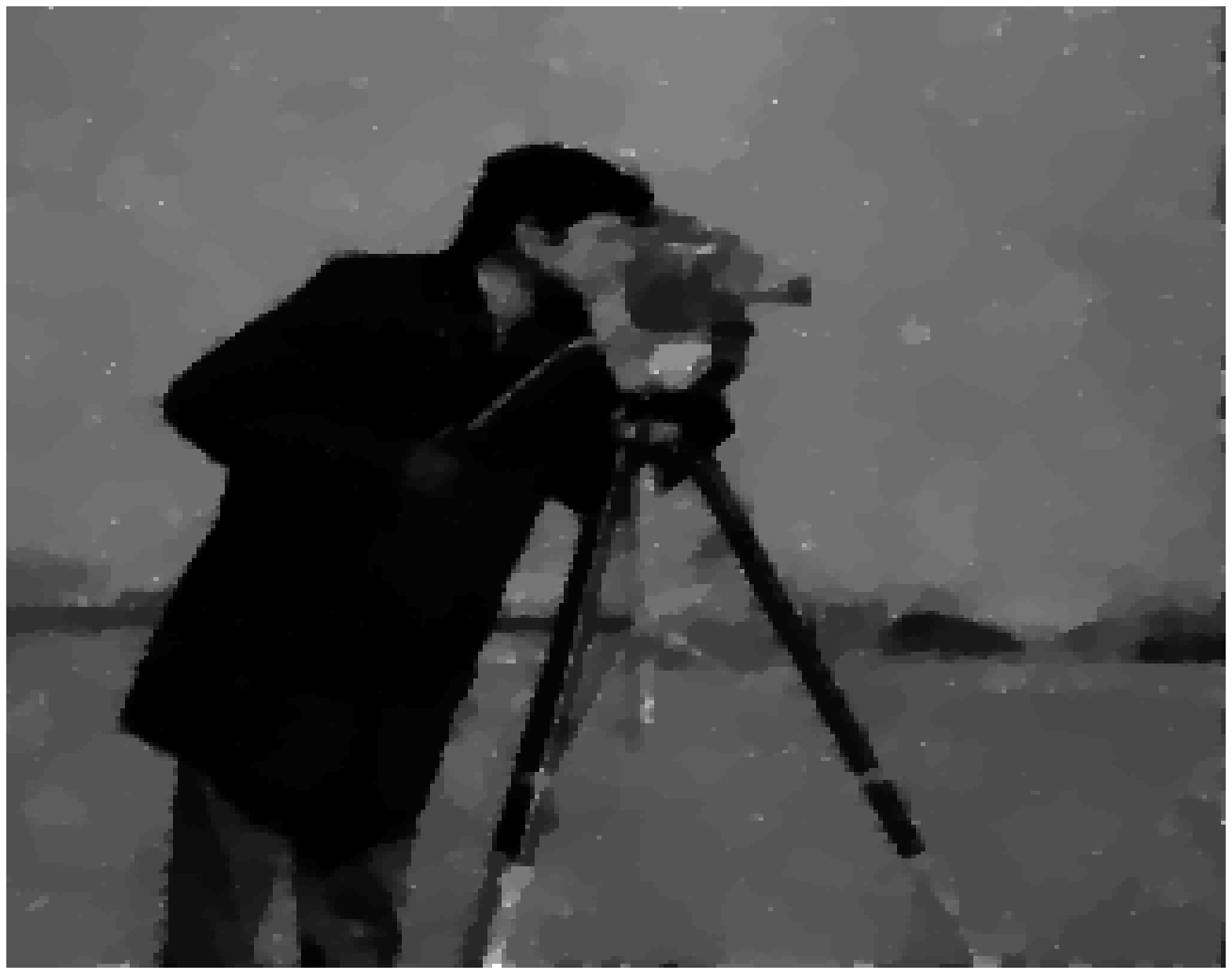}

\vspace{0.1cm}
\includegraphics[width=0.48\columnwidth,height=0.48\columnwidth]{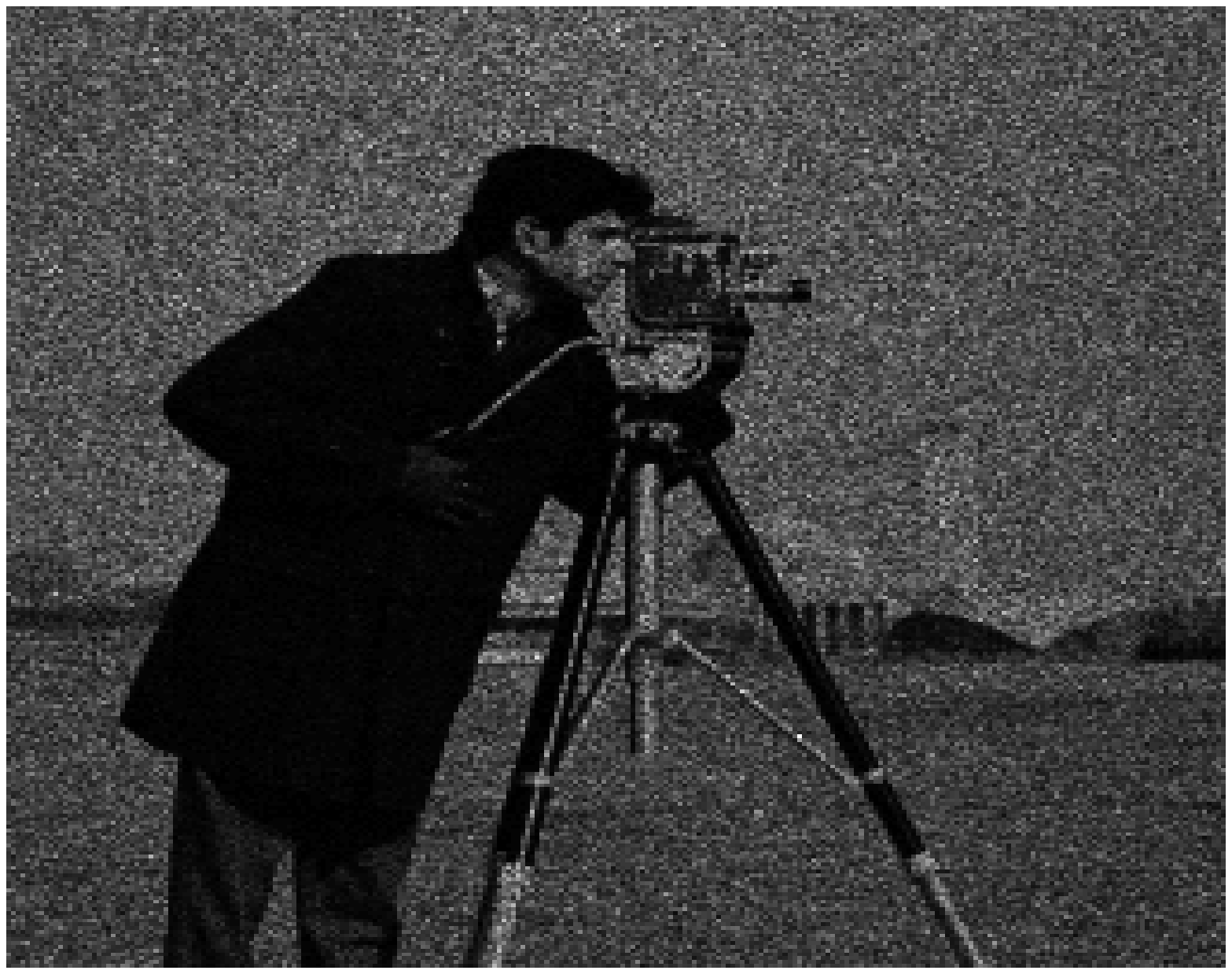}\hspace{.019\columnwidth}\includegraphics[width=0.48\columnwidth,height=0.48\columnwidth]{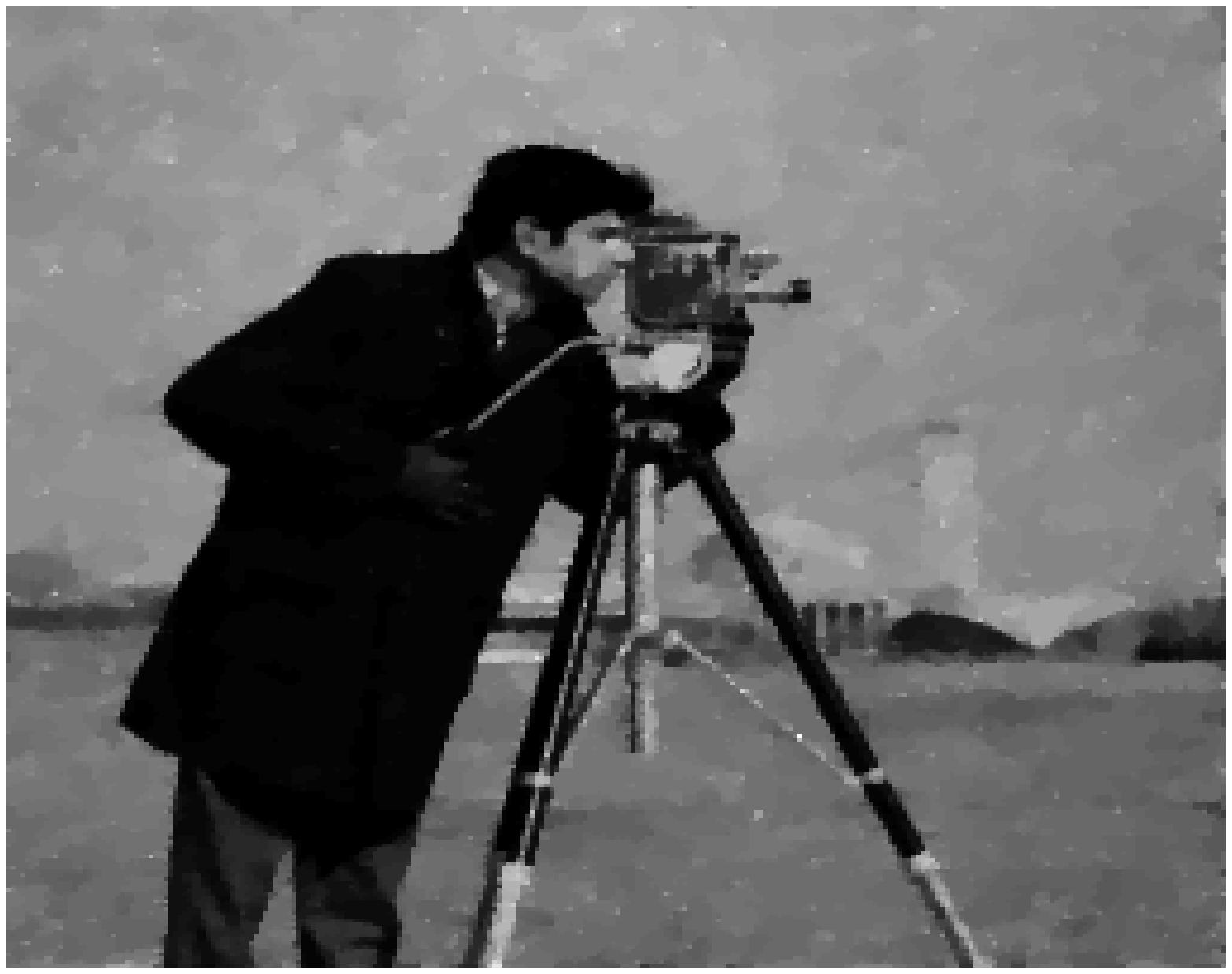}

\vspace{0.1cm}
\includegraphics[width=0.48\columnwidth,height=0.48\columnwidth]{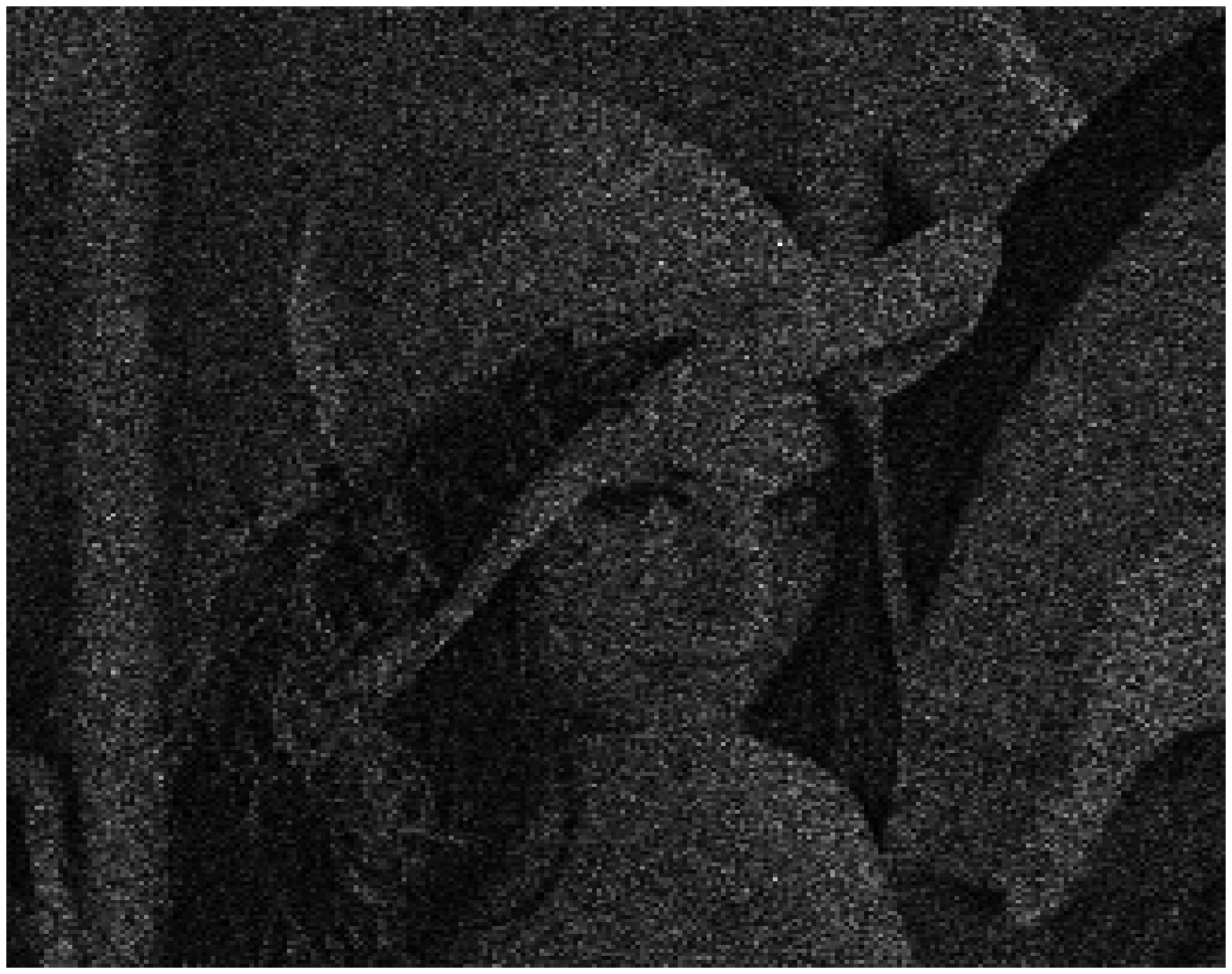}\hspace{.019\columnwidth}\includegraphics[width=0.48\columnwidth,height=0.48\columnwidth]{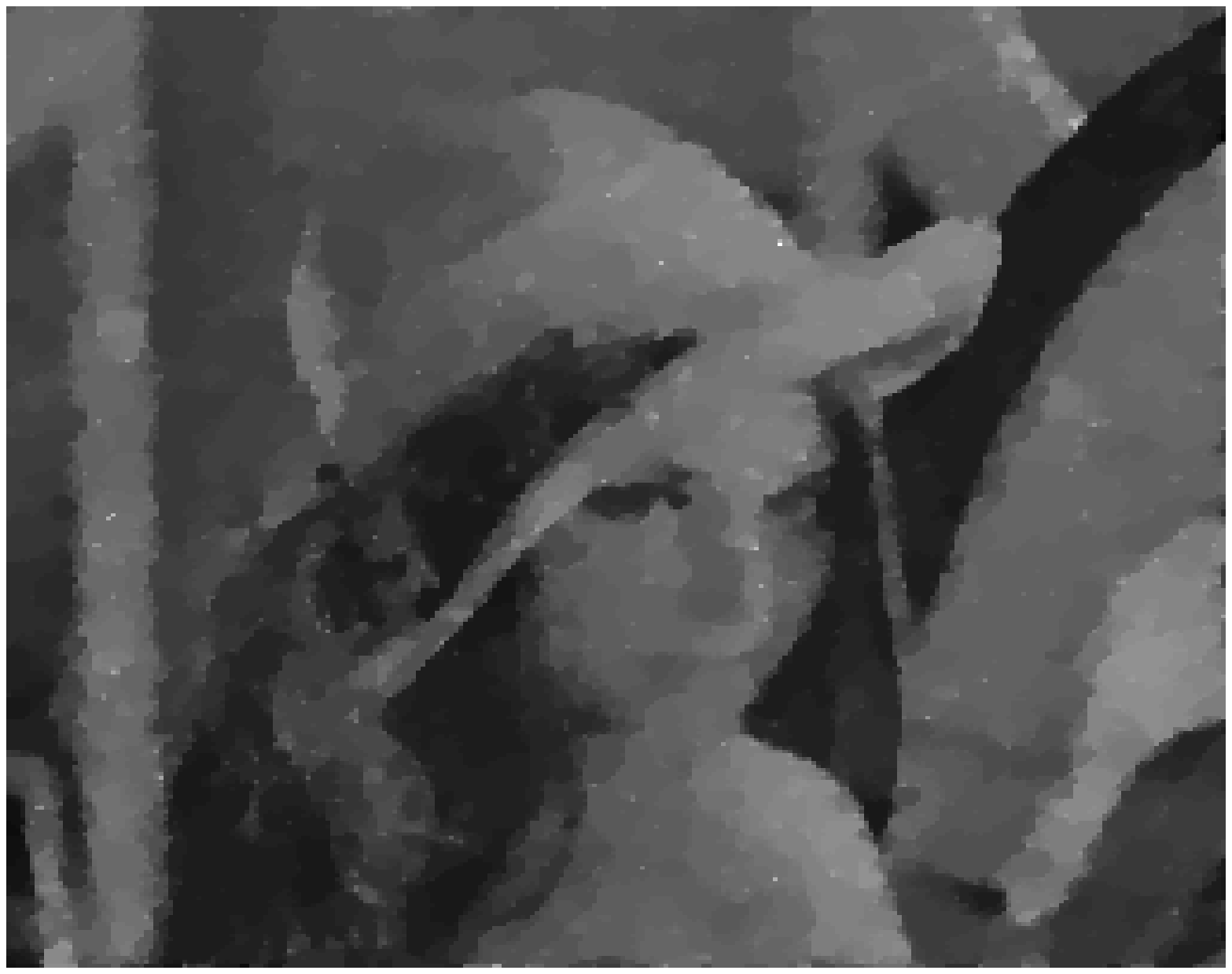}

\vspace{0.1cm}
\includegraphics[width=0.48\columnwidth,height=0.48\columnwidth]{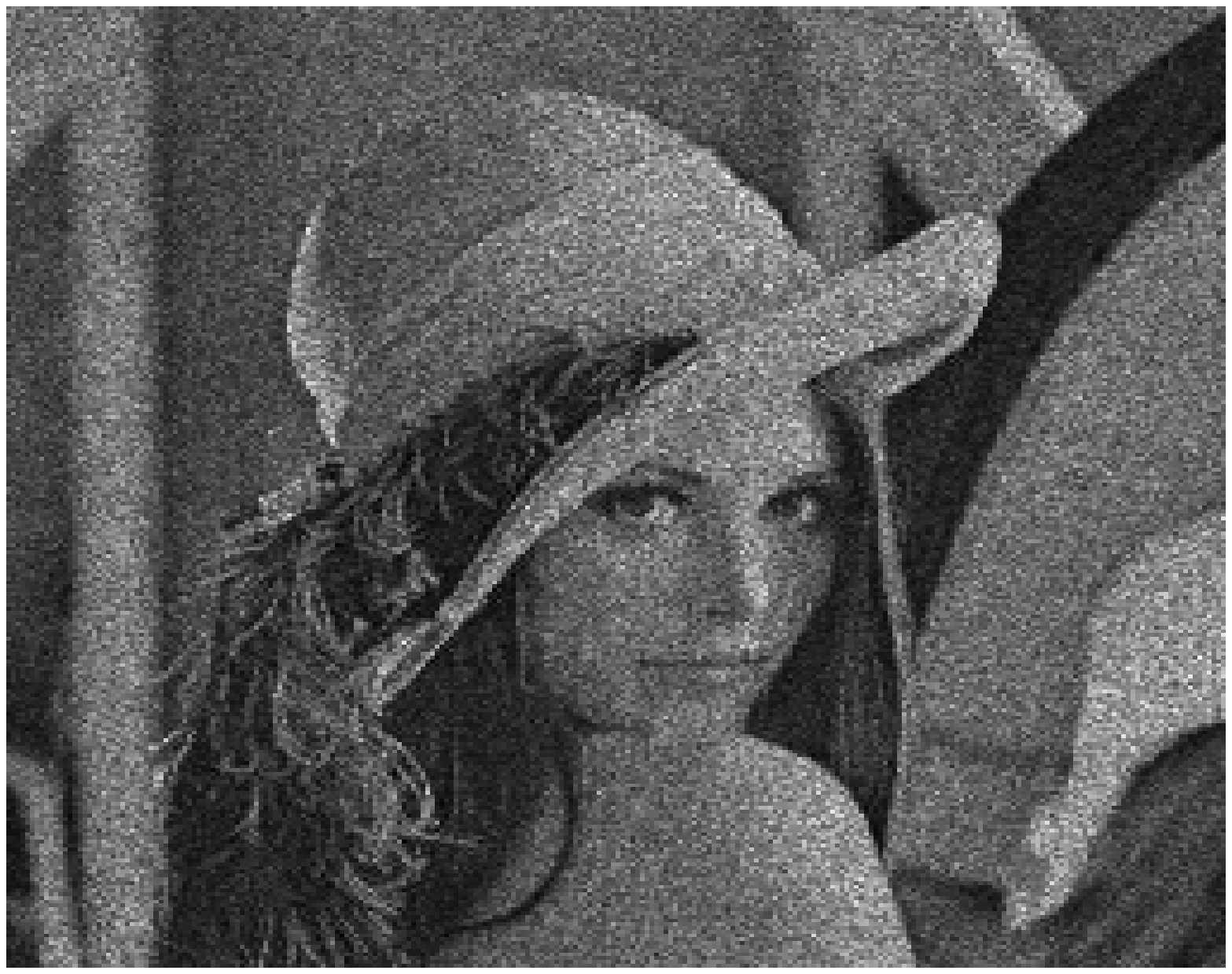}\hspace{.019\columnwidth}\includegraphics[width=0.48\columnwidth,height=0.48\columnwidth]{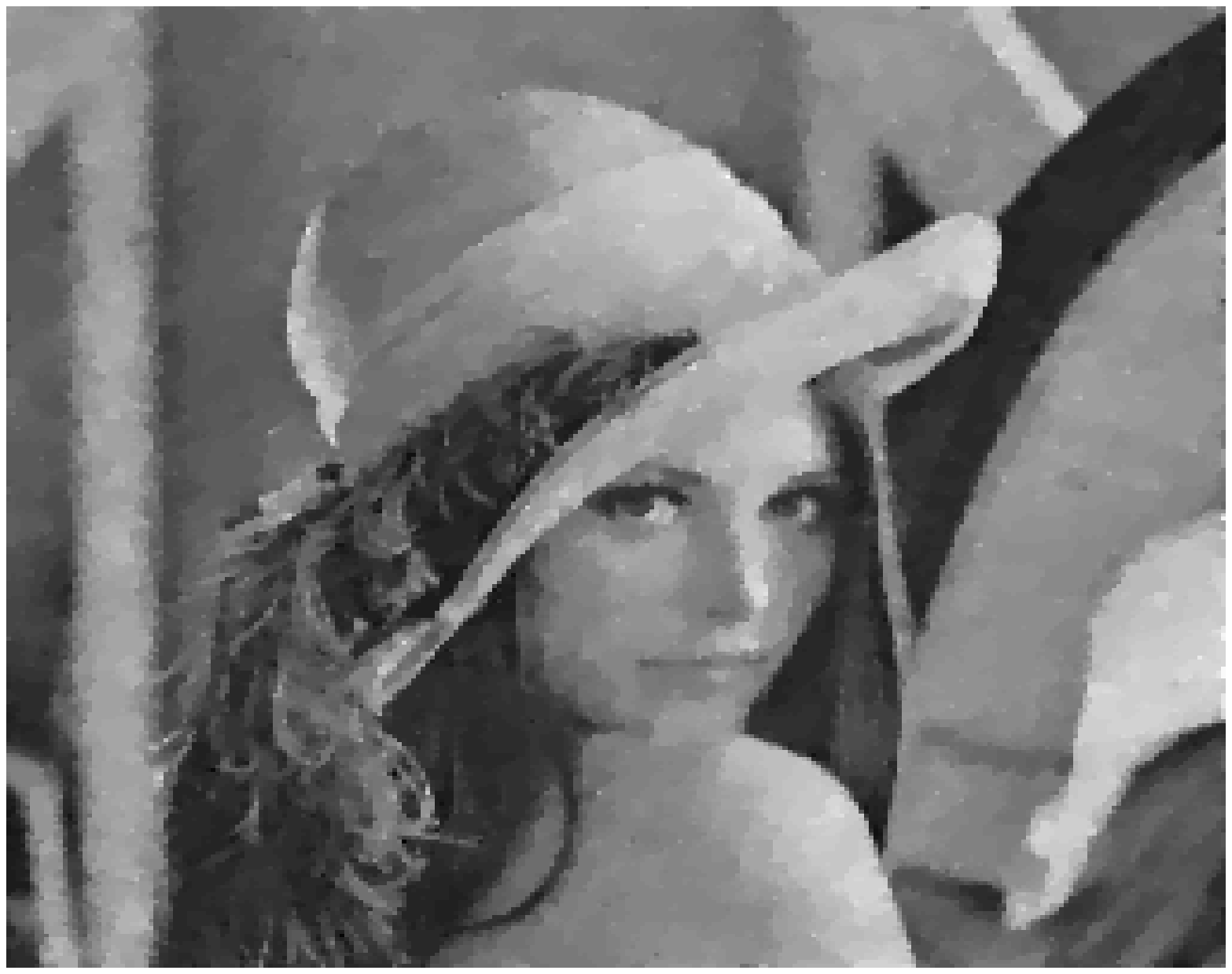}
\caption{Left column: observed noisy images for Experiments 1 to 4 with $M = 3, 13, 5,
33$, respectively. Right column: image estimates.} \label{fig:exps_1_to_4}
\end{figure}

\begin{figure}
\centering
\includegraphics[width=0.48\columnwidth,height=0.48\columnwidth]{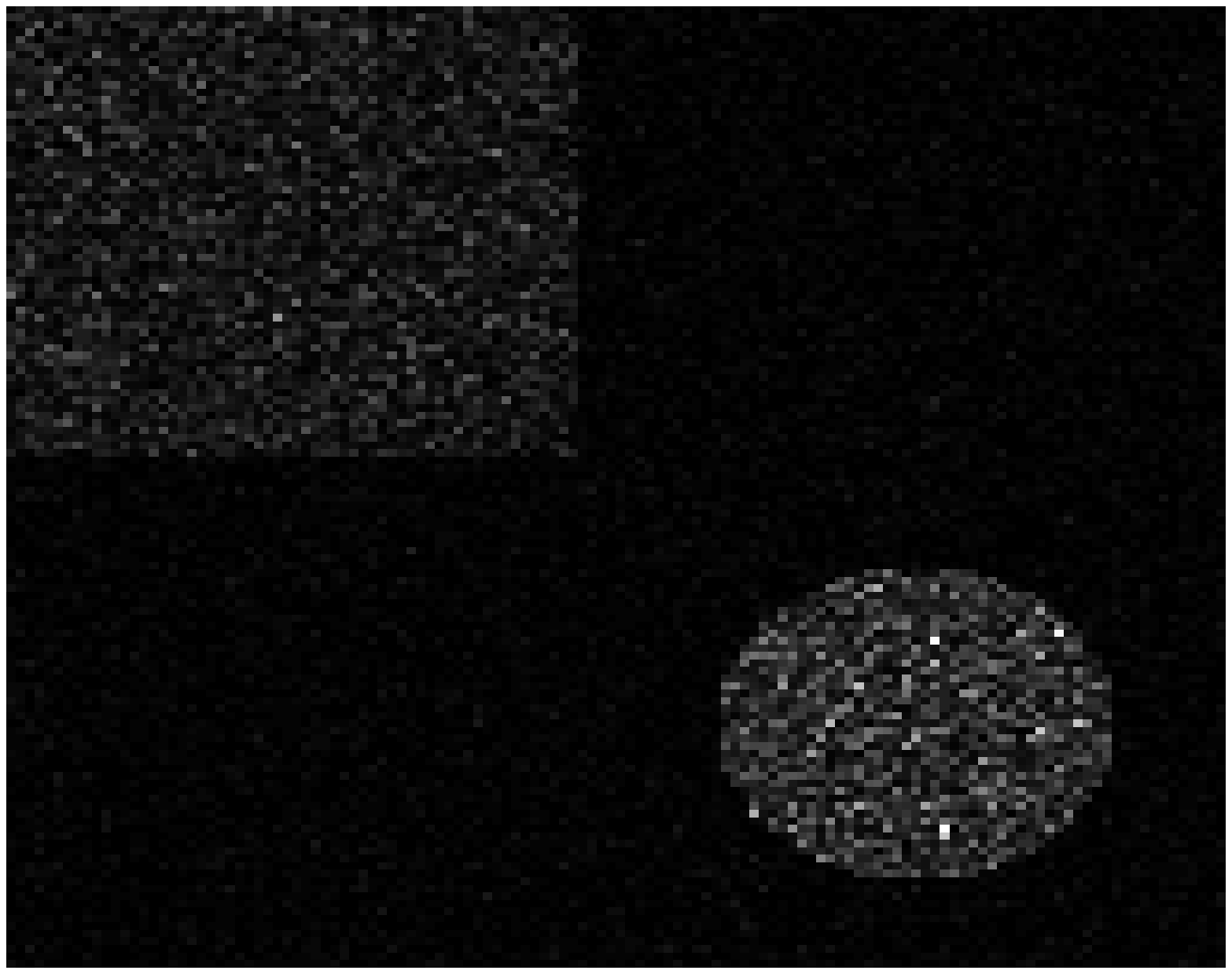}\hspace{.019\columnwidth}\includegraphics[width=0.48\columnwidth,height=0.48\columnwidth]{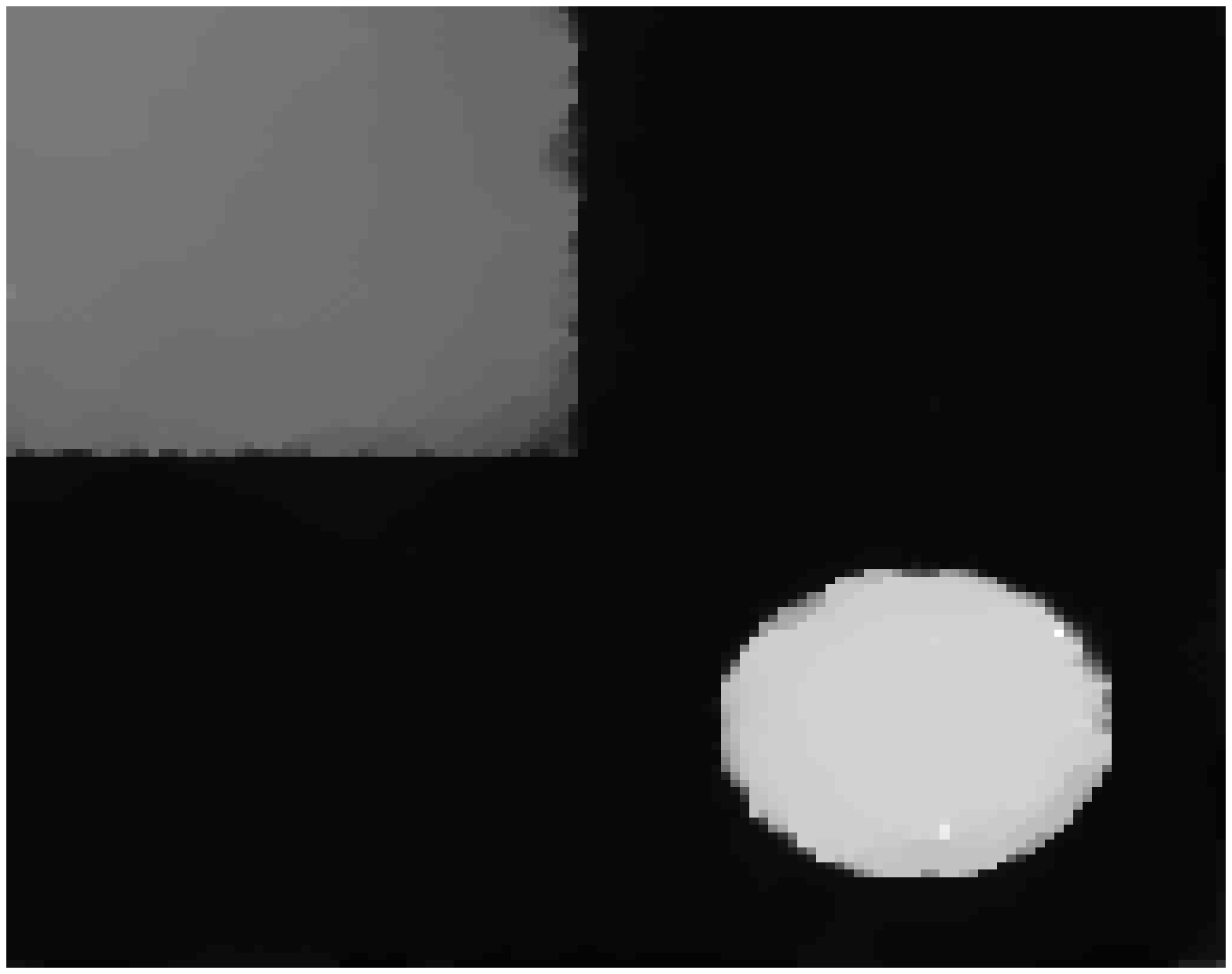}

\vspace{0.1cm}
\includegraphics[width=0.48\columnwidth,height=0.48\columnwidth]{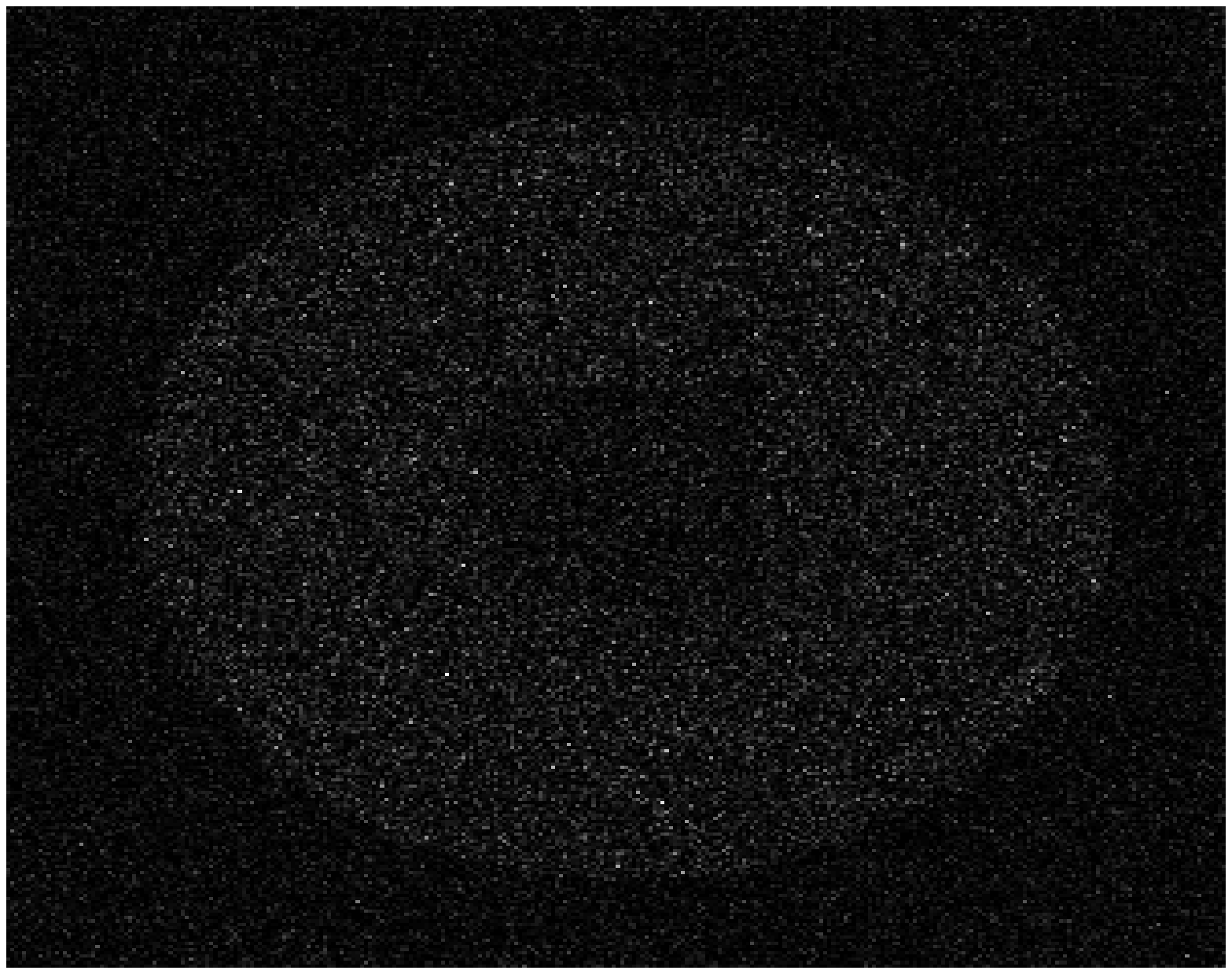}\hspace{.019\columnwidth}\includegraphics[width=0.48\columnwidth,height=0.48\columnwidth]{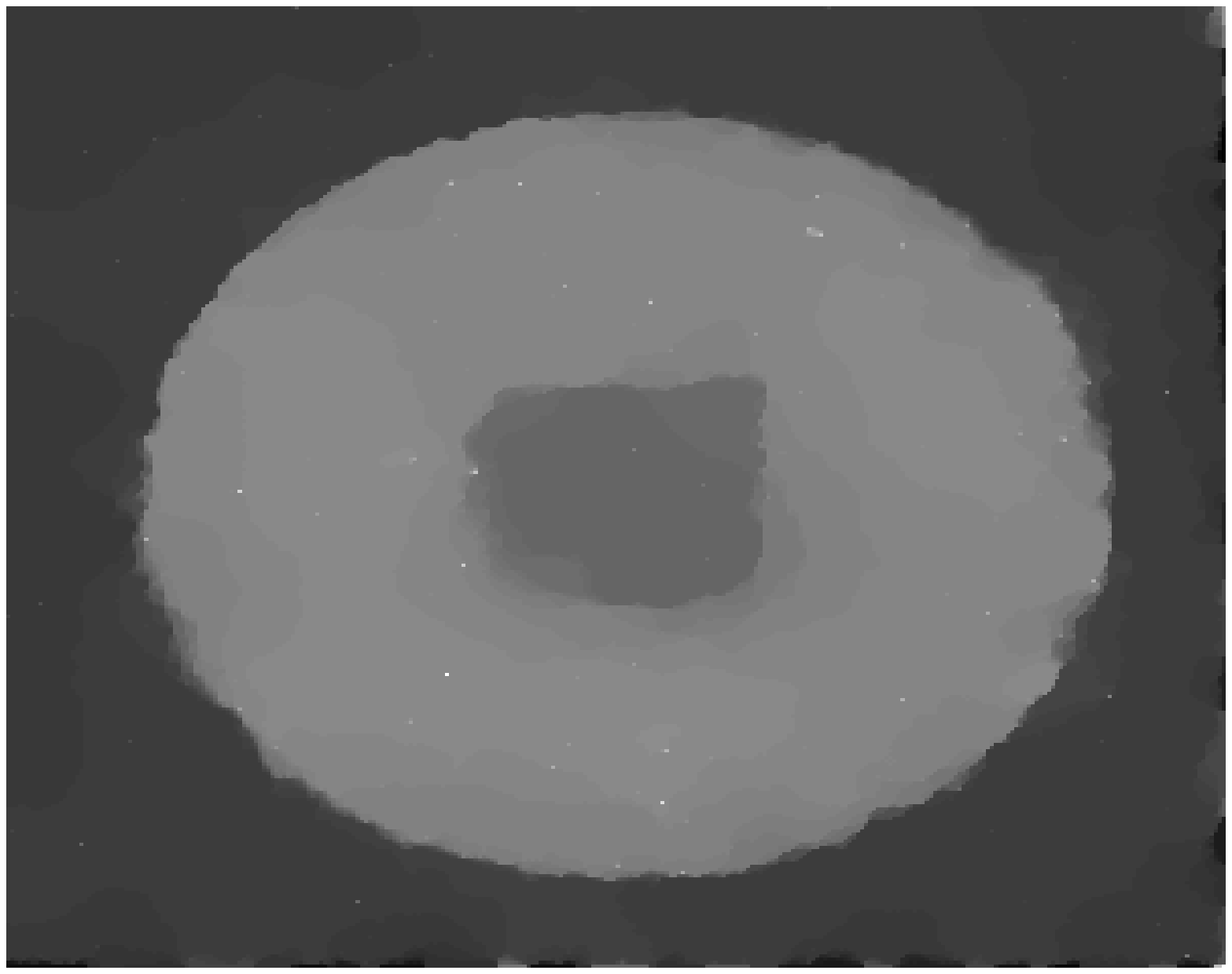}

\vspace{0.1cm}
\includegraphics[width=0.48\columnwidth,height=0.48\columnwidth]{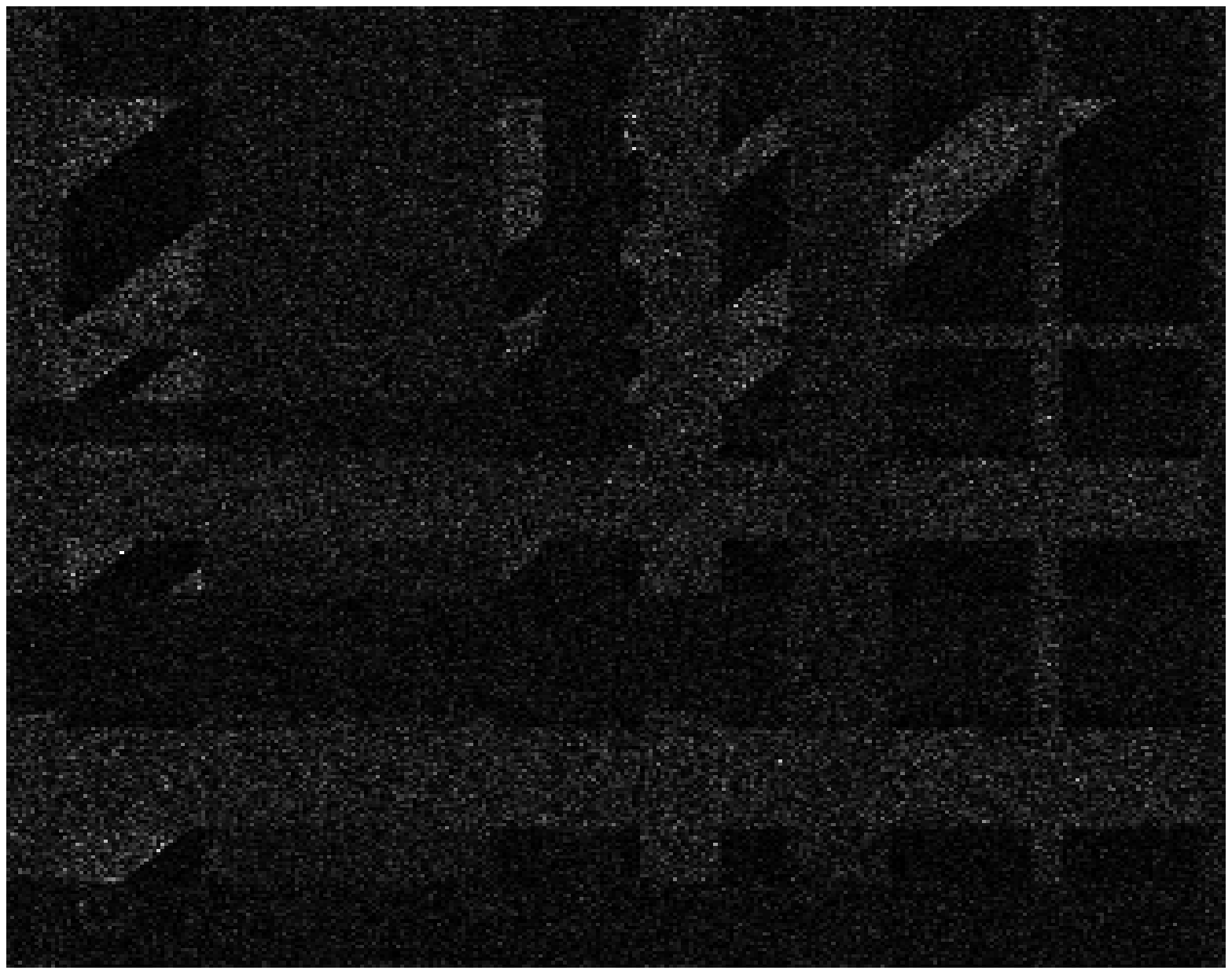}\hspace{.019\columnwidth}\includegraphics[width=0.48\columnwidth,height=0.48\columnwidth]{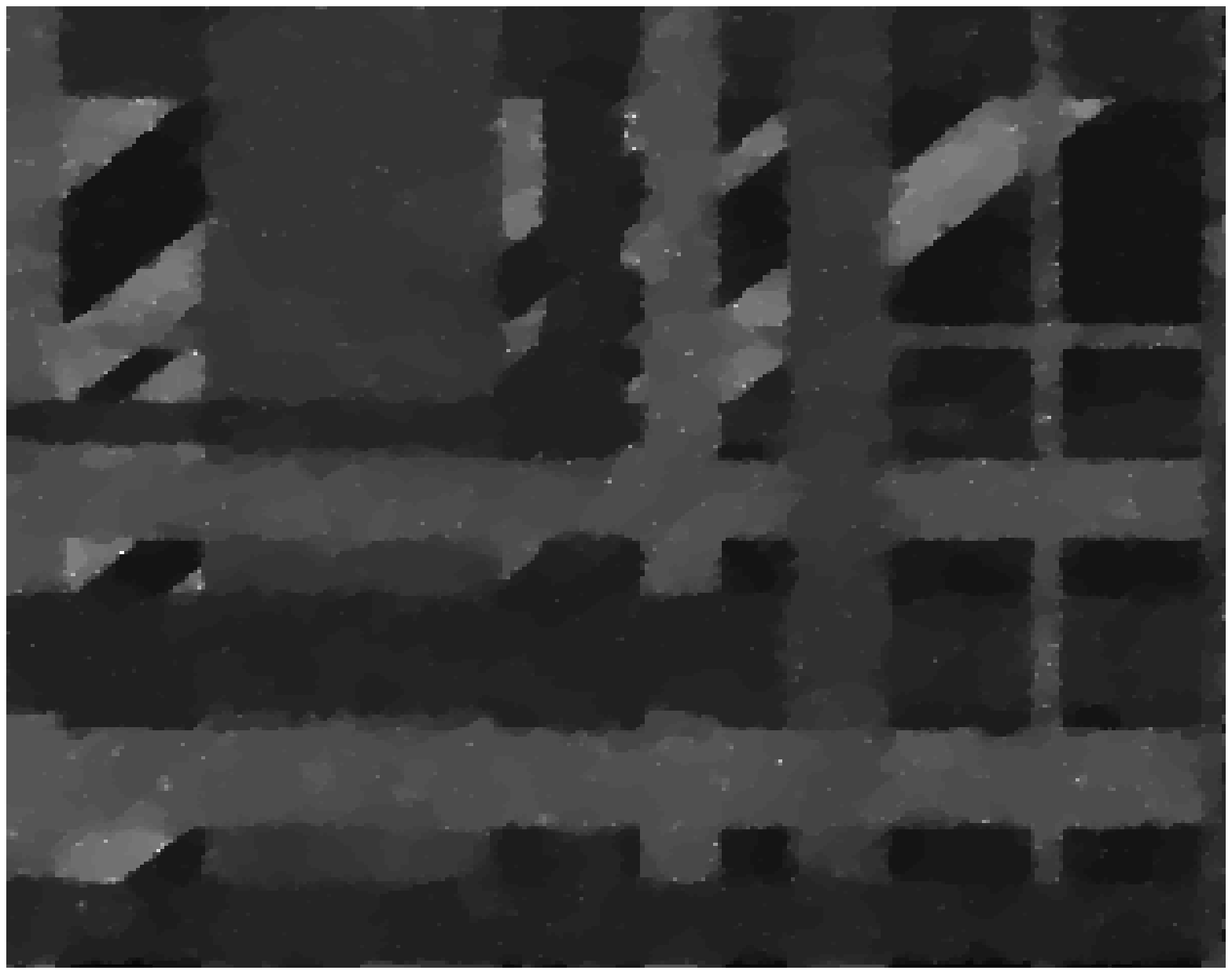}

\caption{Left column: observed noisy images for Experiments 5 to 7 with $M = 2, 1, 2$, respectively. Right column: image estimates.} \label{fig:exps_5_to_7}
\end{figure}

\begin{figure}

\includegraphics[width=0.48\columnwidth,height=0.48\columnwidth]{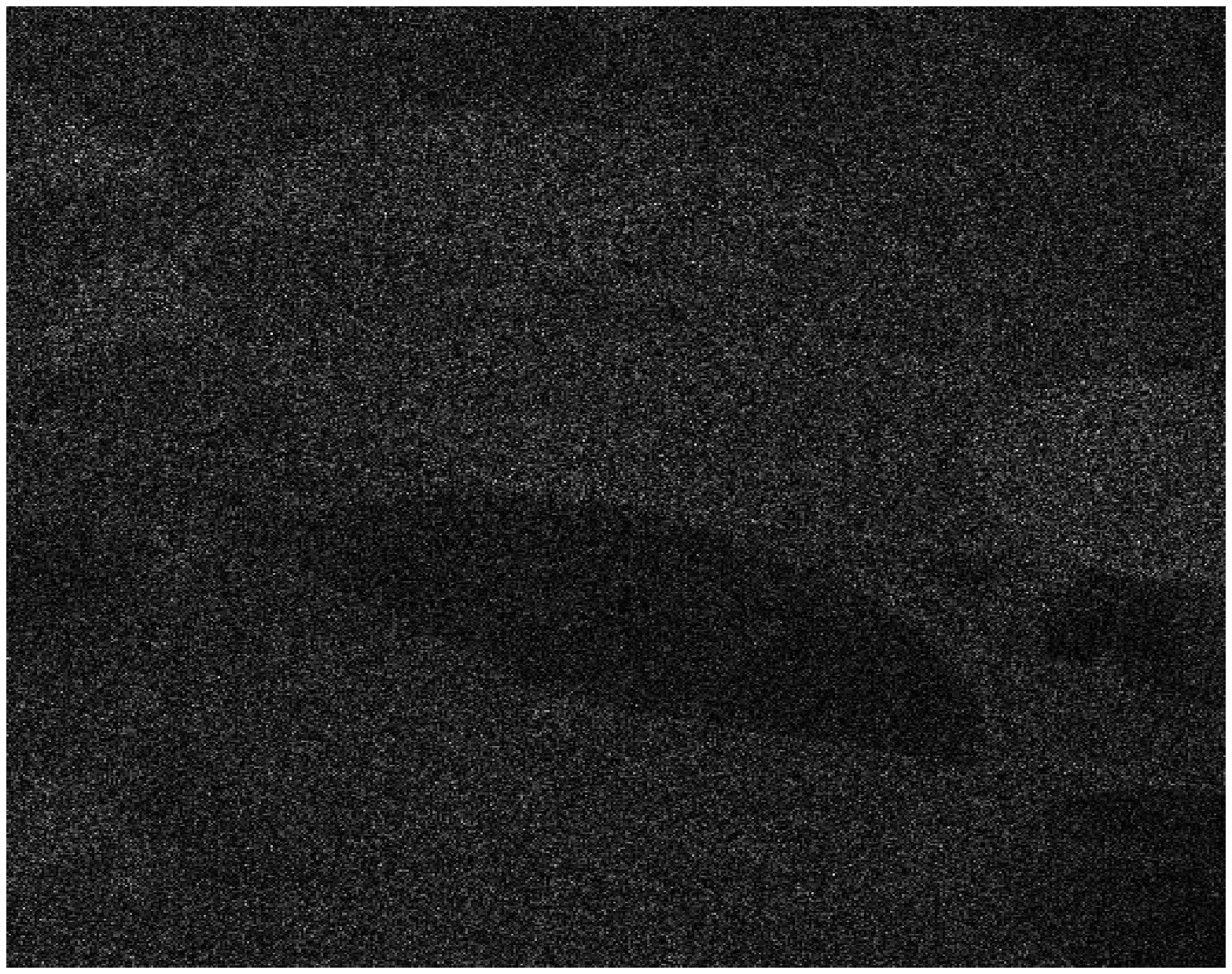}\hspace{.019\columnwidth}\includegraphics[width=0.48\columnwidth,height=0.48\columnwidth]{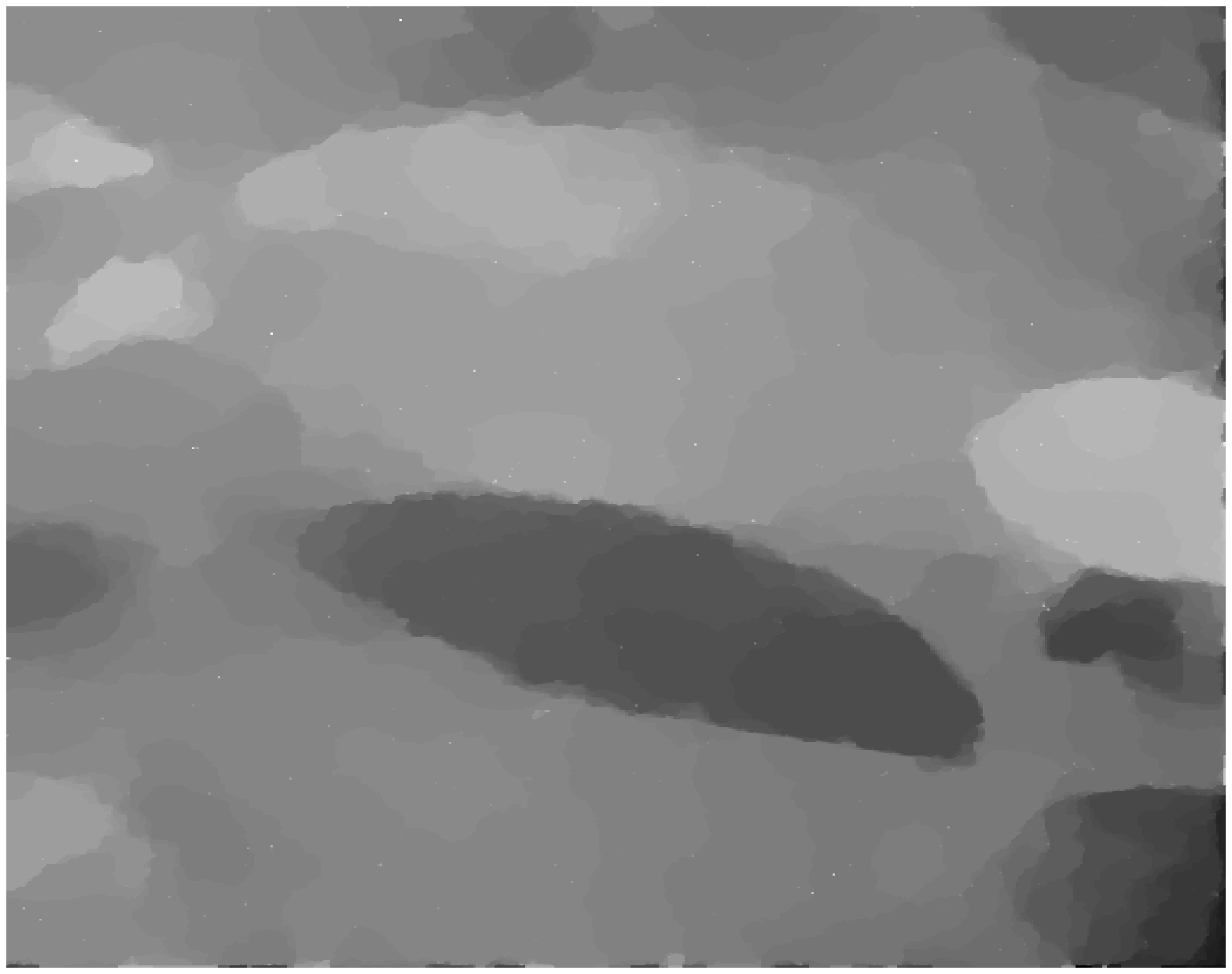}

\vspace{0.1cm}
\includegraphics[width=0.48\columnwidth,height=0.48\columnwidth]{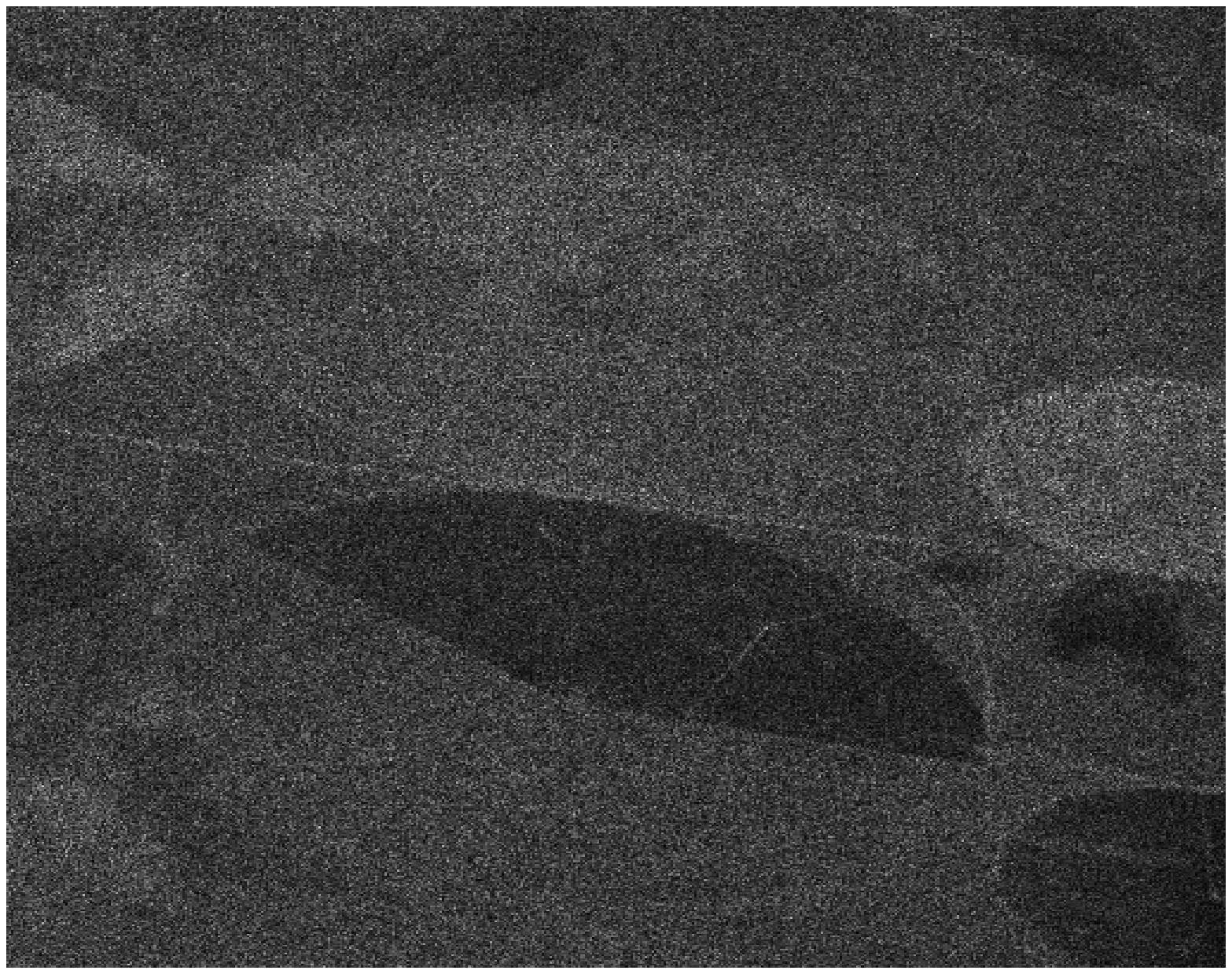}\hspace{.019\columnwidth}\includegraphics[width=0.48\columnwidth,height=0.48\columnwidth]{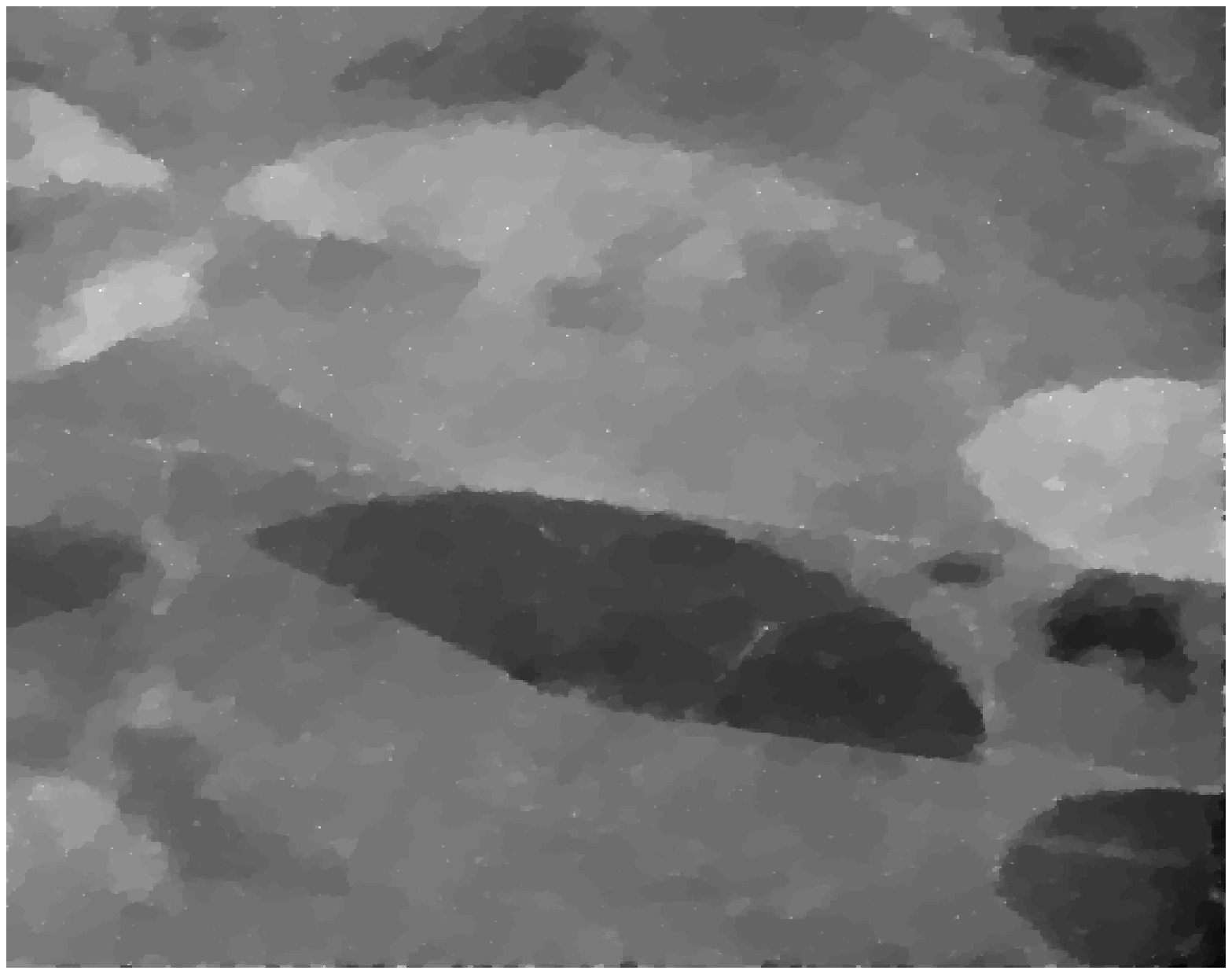}

\vspace{0.1cm}
\includegraphics[width=0.48\columnwidth,height=0.48\columnwidth]{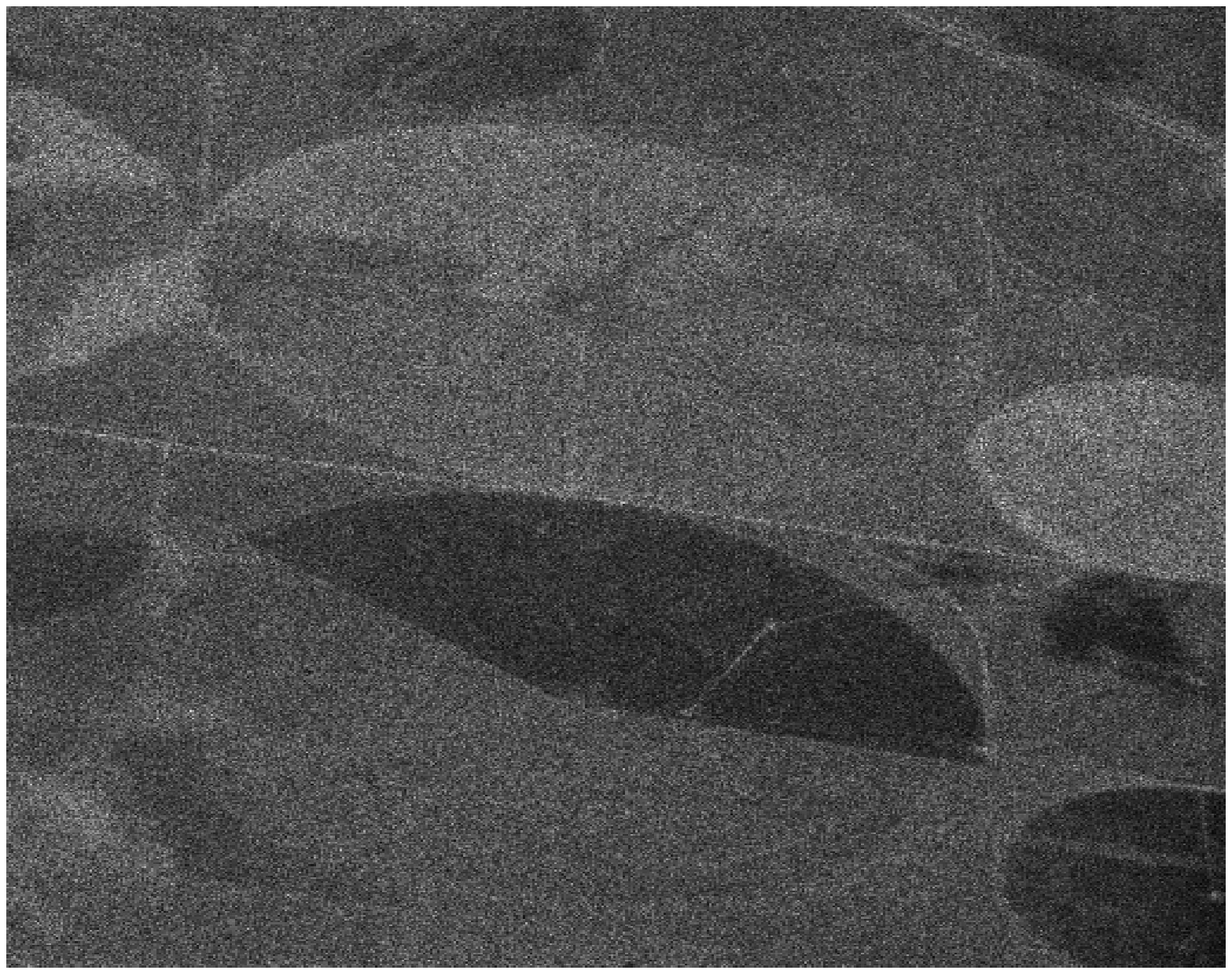}\hspace{.019\columnwidth}\includegraphics[width=0.48\columnwidth,height=0.48\columnwidth]{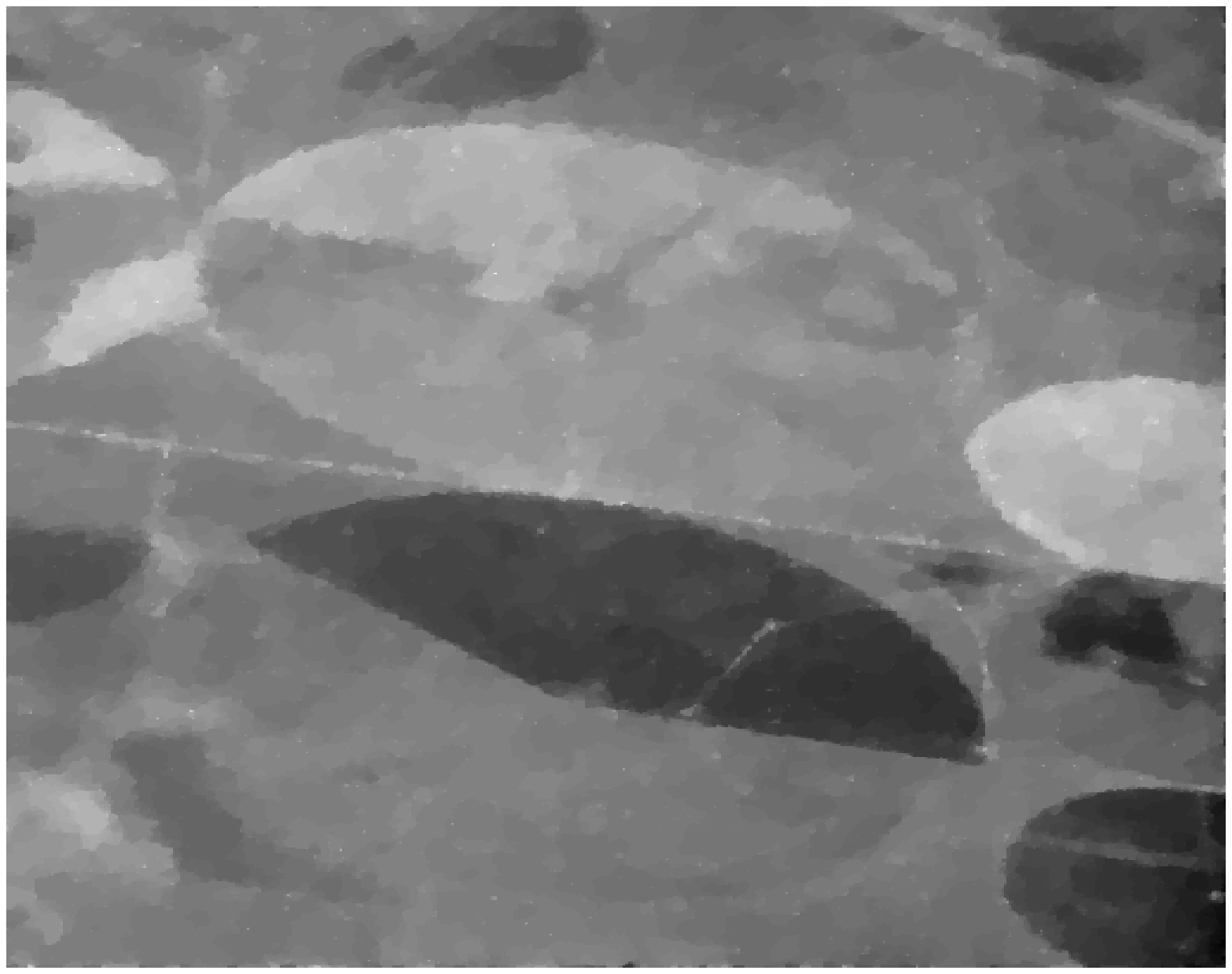}

\caption{Left column: observed noisy images for Experiments 8 to 10. Right column: image estimates. Note: for better
visualization, all the images underwent the nonlinear transformation $(\cdot)^{0.7}$ prior to being displayed.} \label{fig:exps_8_to_10}
\end{figure}

\begin{figure}

\includegraphics[width=0.48\columnwidth,height=0.48\columnwidth]{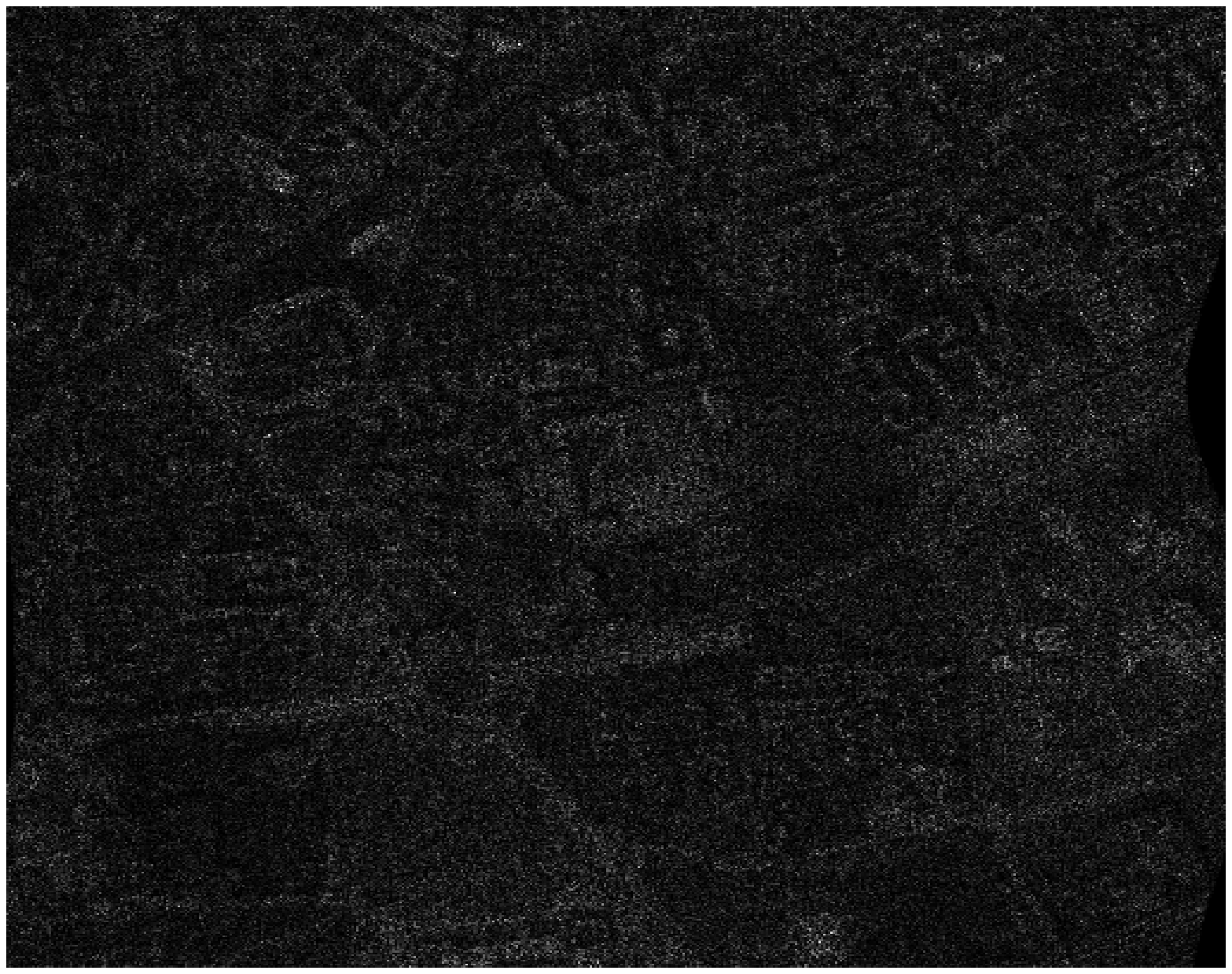}\hspace{.019\columnwidth}\includegraphics[width=0.48\columnwidth,height=0.48\columnwidth]{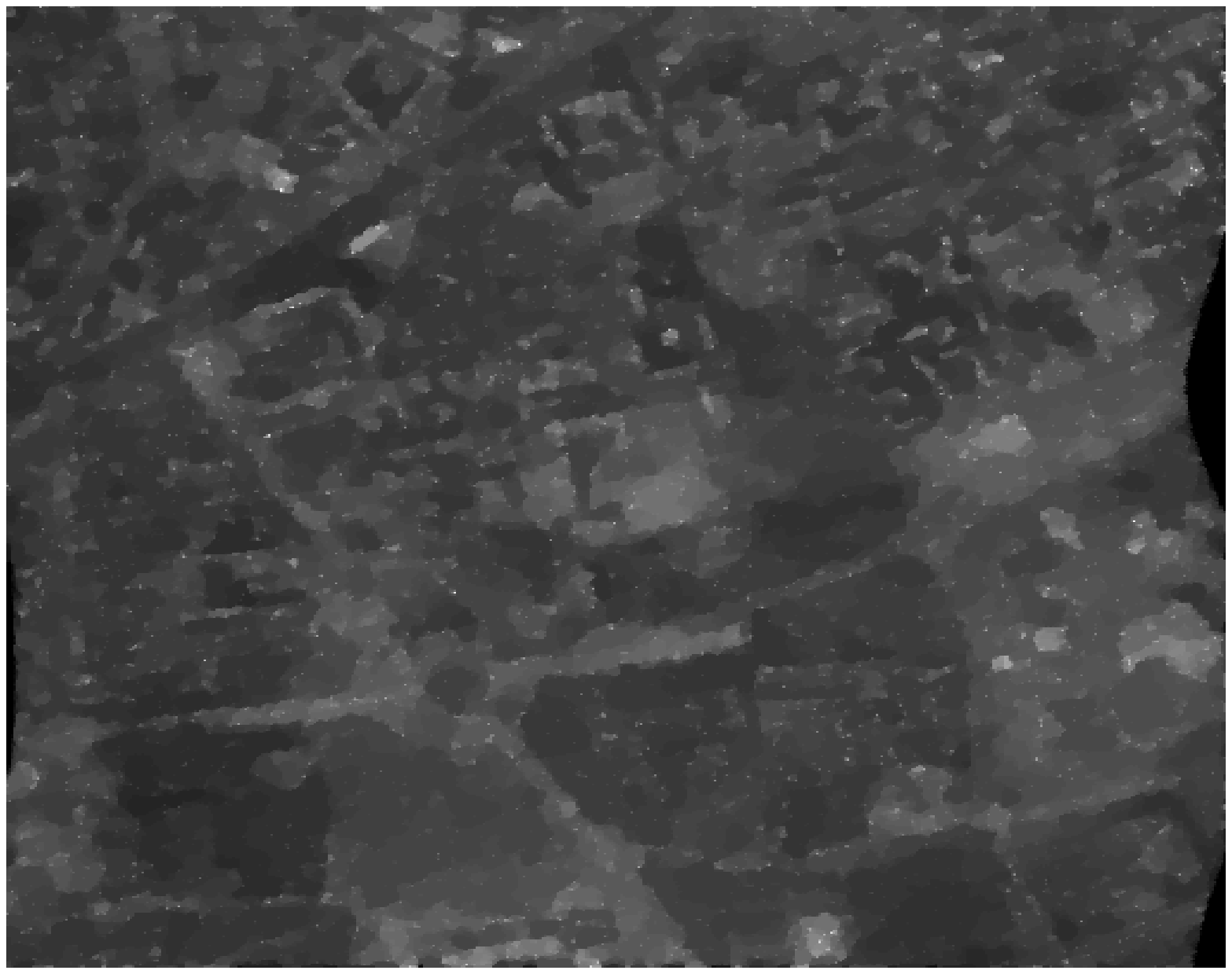}

\vspace{0.1cm}
\includegraphics[width=0.48\columnwidth,height=0.48\columnwidth]{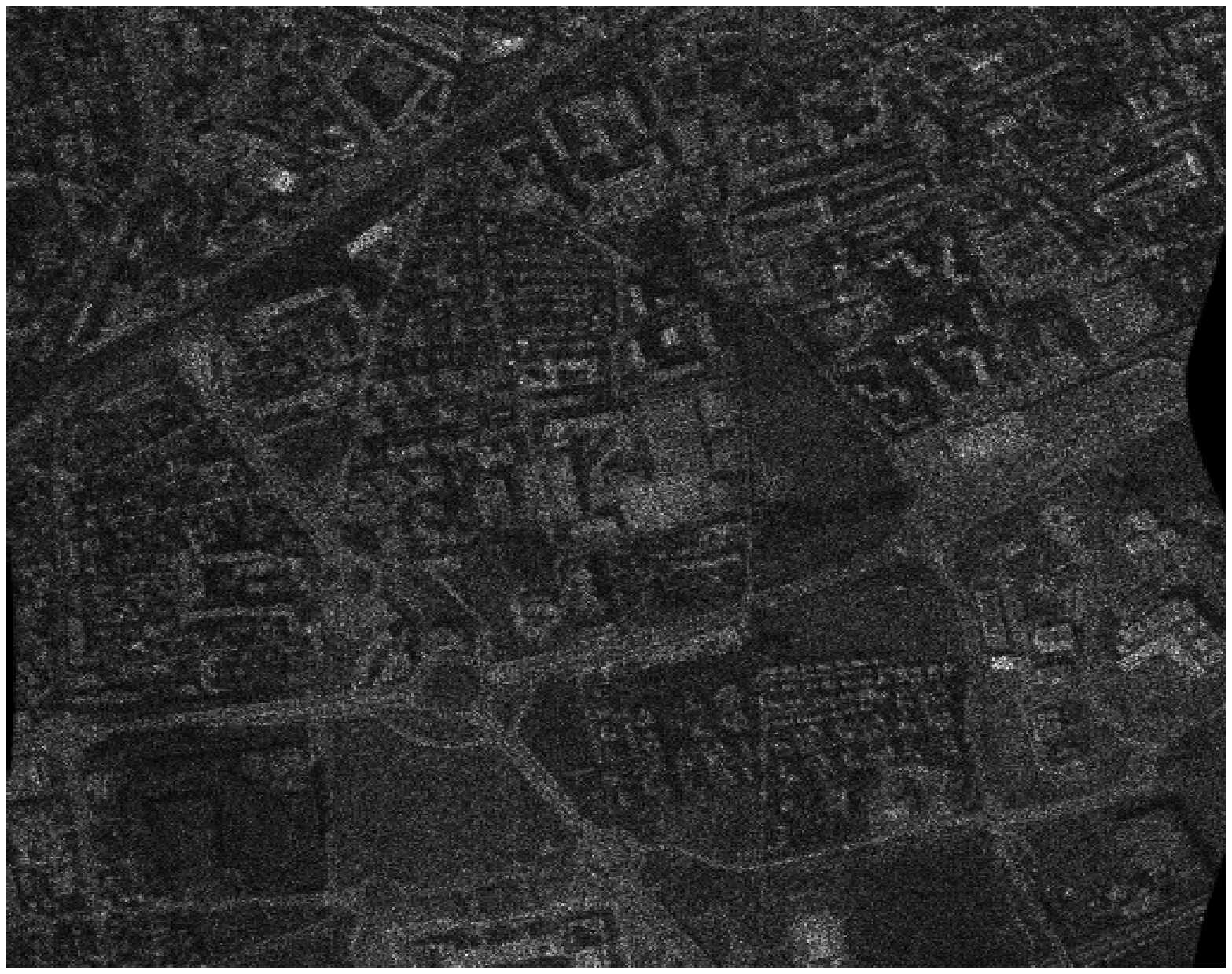}\hspace{.019\columnwidth}\includegraphics[width=0.48\columnwidth,height=0.48\columnwidth]{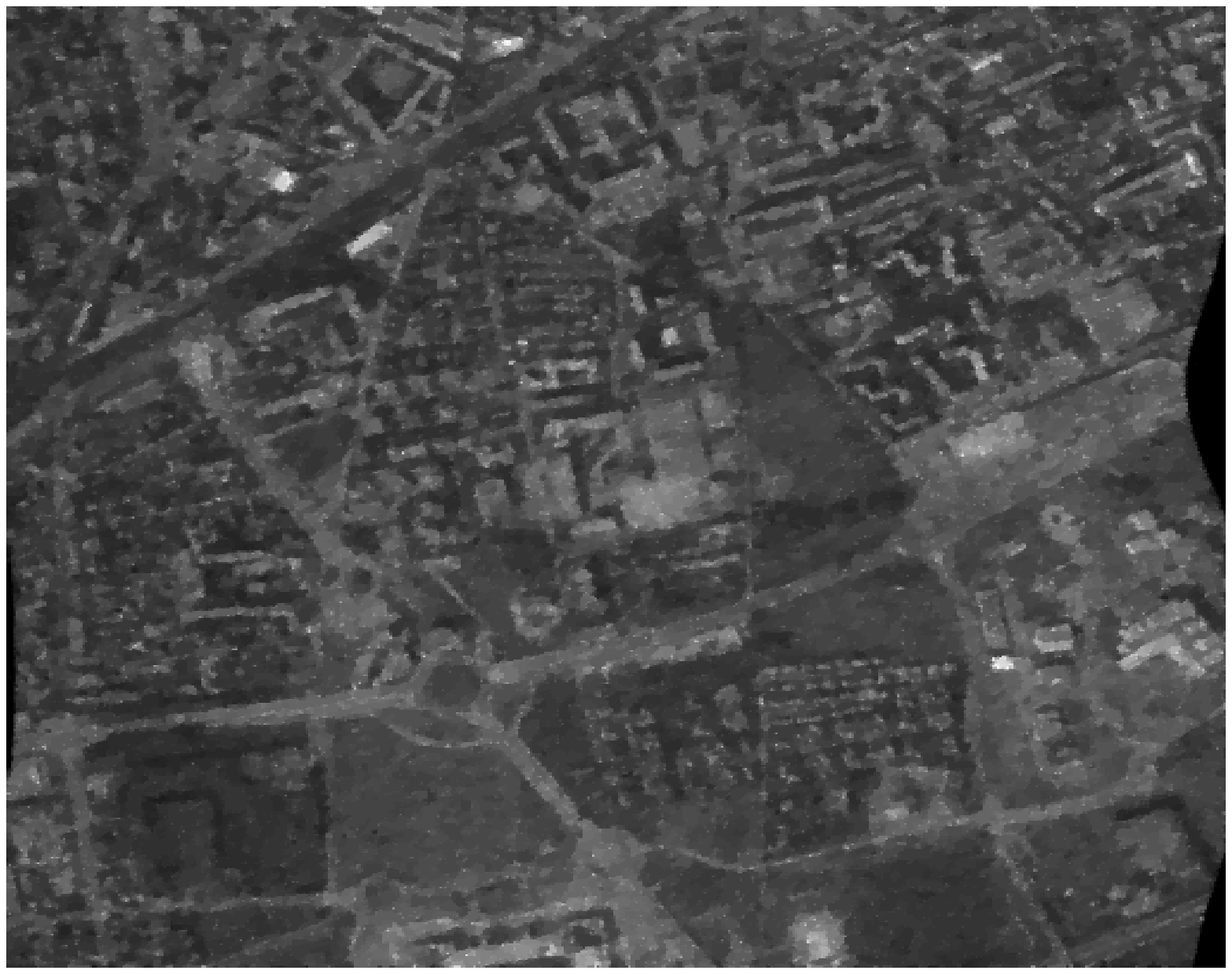}

\vspace{0.1cm}
\includegraphics[width=0.48\columnwidth,height=0.48\columnwidth]{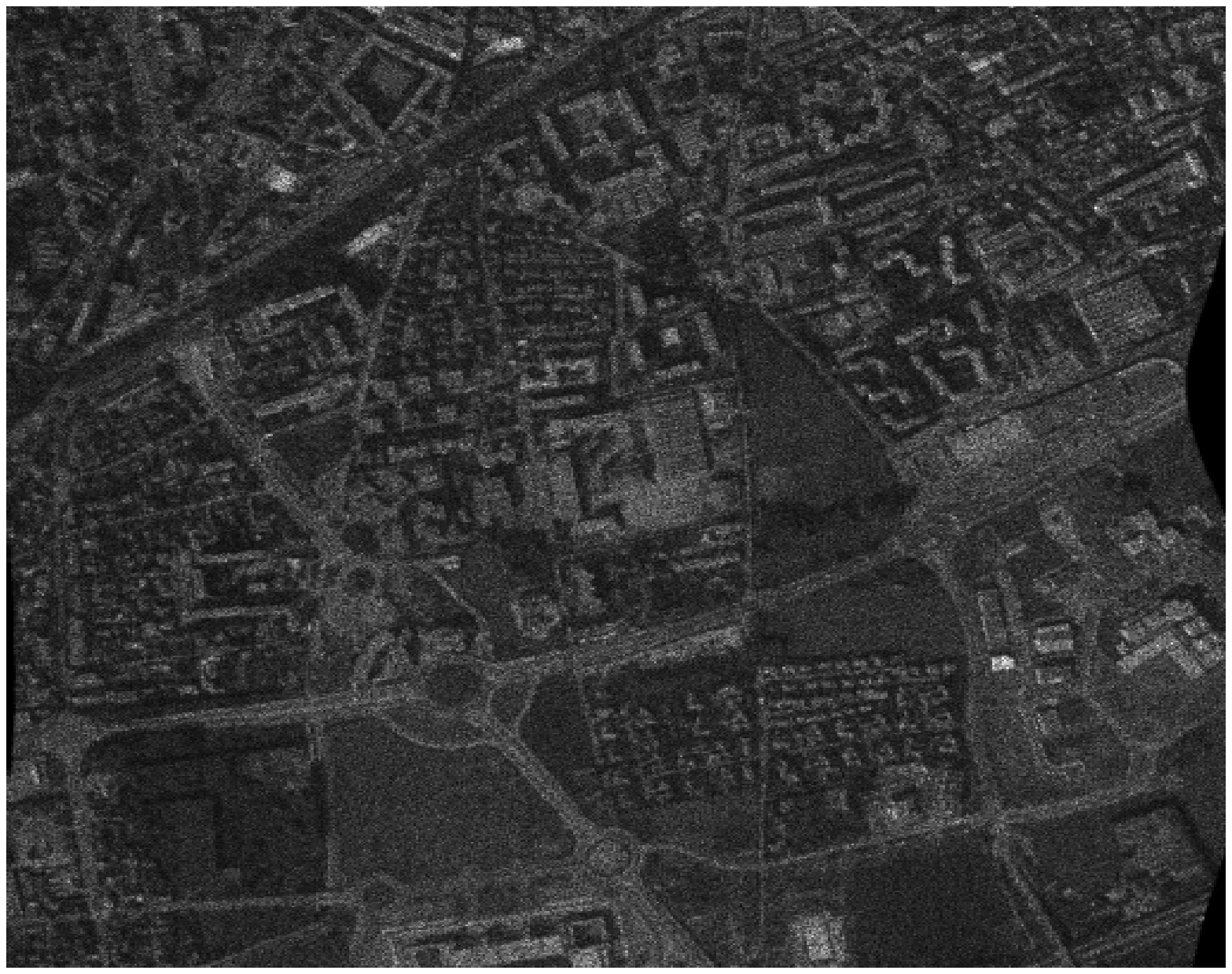}\hspace{.019\columnwidth}\includegraphics[width=0.48\columnwidth,height=0.48\columnwidth]{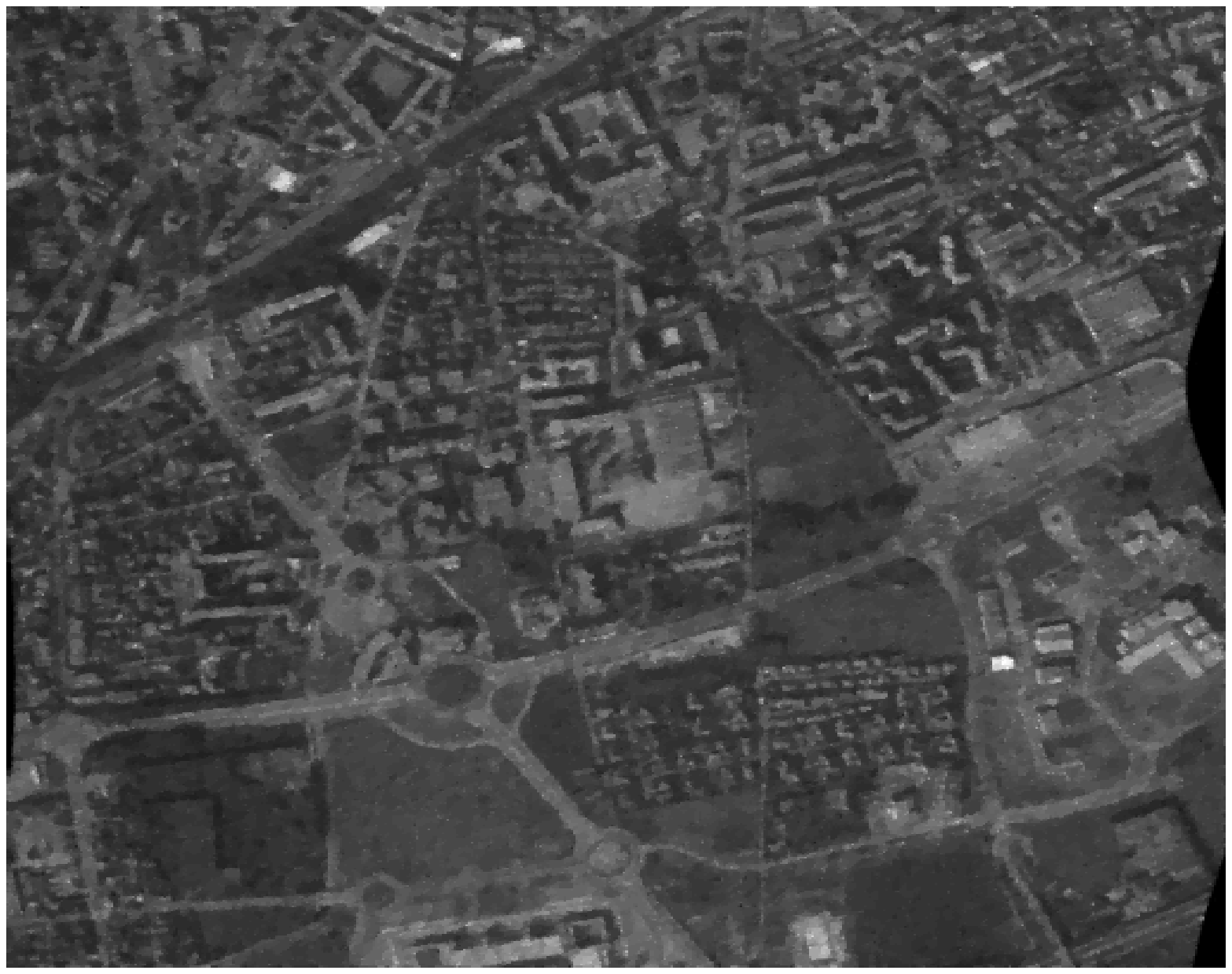}

\caption{Left column: observed noisy images for Experiments 11 to 13. Right column: image estimates.
Note: for better visualization, all the images underwent the nonlinear transformation $(\cdot)^{0.7}$ prior to being displayed.} \label{fig:exps_11_to_13}
\end{figure}

\begin{figure}
\centerline{\includegraphics[width=5.8cm]{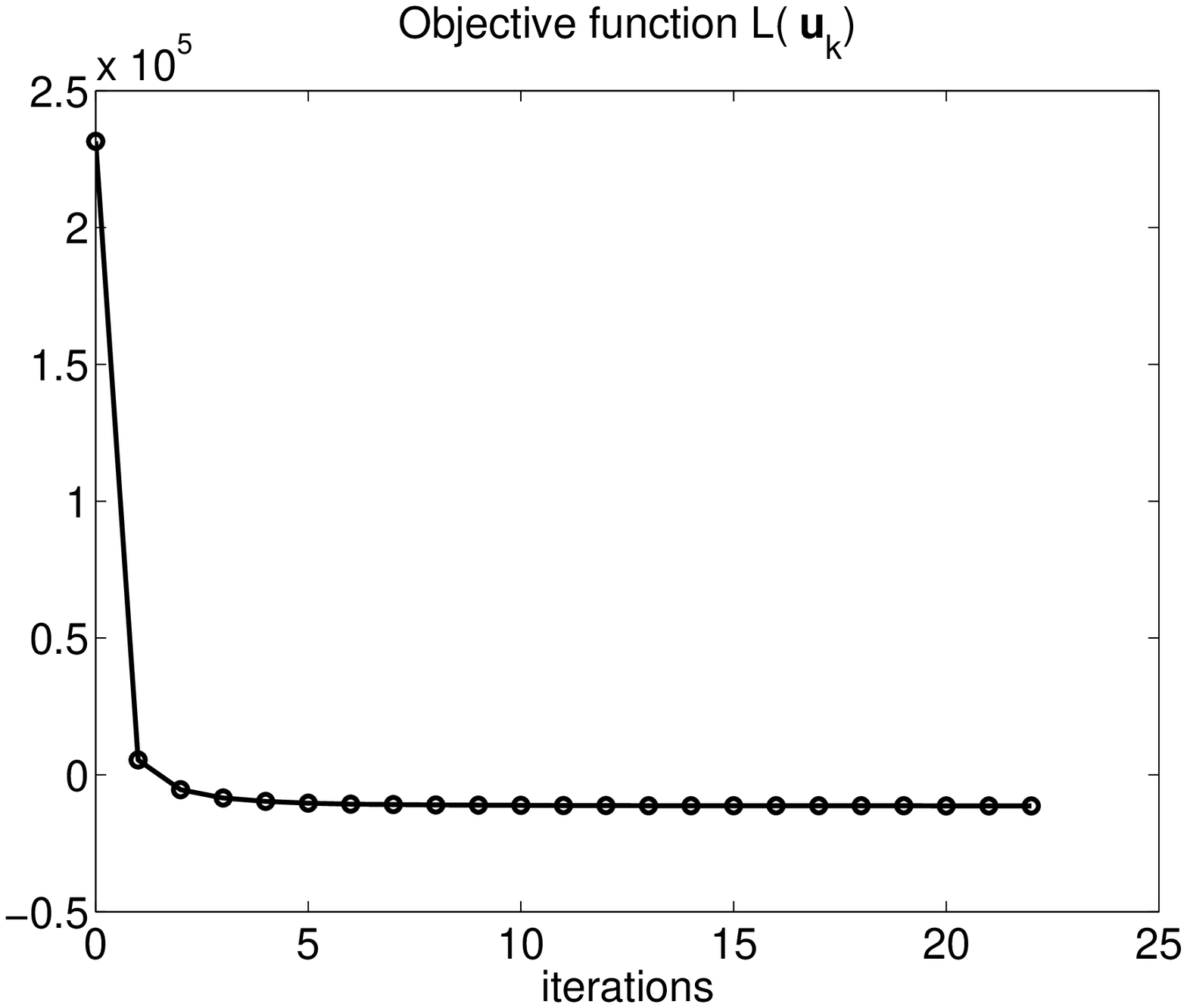}}

\centerline{\includegraphics[width=5.8cm]{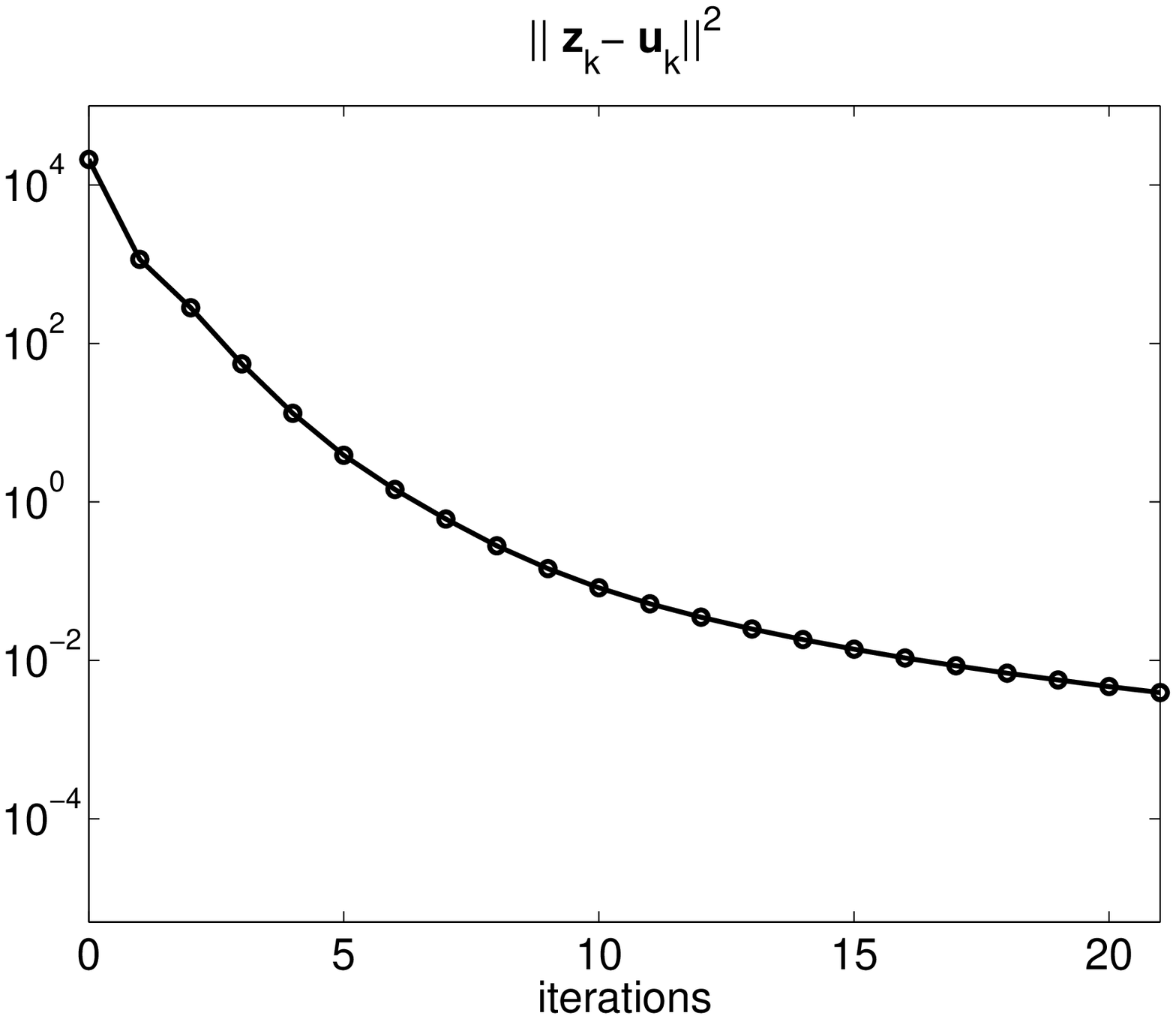}}
\caption{Evolution of the objective function $L({\bf u}_k)$ and of the constraint
function $\|{\bf z}_k - {\bf u}_k\|_2^2$, along the iterations of the algorithm, for
Experiment 1. } \label{fig:plots}
\end{figure}

Figures \ref{fig:exps_1_to_4}, \ref{fig:exps_5_to_7},  \ref{fig:exps_8_to_10}, and
\ref{fig:exps_11_to_13} show the noisy and restored images, for the  Experiments 1 to 4, 5 to 7,  8 to 10,
and 11 to 13, respectively. Finally, Fig.~\ref{fig:plots} plots the evolution of the
objective function $L({\bf u}_k)$ and of the constraint function $\|{\bf z}_k - {\bf
u}_k\|_2^2$ along the iterations, for the Experiment 1 (Cameraman image and $M=3$).
notice the decrease of approximately 7 orders of magnitude of  $\|{\bf z}_k - {\bf
u}_k\|_2^2$ along the 21 MIDAL  iterations, showing that, for all practical purposes, the
constraint (\ref{eq:z_u_c}) is satisfied.

\section{Concluding Remarks}
\label{sec:conclusion} We have proposed a new approach to solve the optimization
problem resulting from variational (equivalently MAP) estimation of images
observed under multiplicative noise models. Although the proposed formulation
and algorithm can be used with other priors (namely, frame-based), here we
have focused on total-variation regularization. Our approach is based on
two building blocks: (1) the original unconstrained optimization problem
was first transformed into an equivalent
constrained problem, via a variable splitting procedure; (2) this constrained problem was
then addressed using an augmented Lagrangian method, more specifically, the alternating
direction method of multipliers (ADMM). We have shown that the conditions for the convergence of
ADMM are satisfied.

Multiplicative noise removal (equivalently reflectance estimation) was formulated
with respect to the logarithm of the reflectance, as proposed by some other authors.
 As a consequence, the multiplicative
noise was converted into additive noise yielding a strictly convex data term ({\em i.e.},
negative of the log-likelihood function), which was not the case with the original multiplicative
noise model. A consequence of this strict convexity, together with the convexity of
the total variation regularizer, was that the solution of the variational problem
(the denoised image) is unique and the resulting algorithm, termed MIDAL
({\it multiplicative image denoising by augmented Lagrangian}), is guaranteed to converge.

MIDAL  is very simple and, in the experiments herein reported, exhibited state-of-the-art
estimation performance and speed. For example, compared with the hybrid method
in \cite{DurandFadiliN08}, which combines curvelet-based and total-variation
regularization, MIDAL yields comparable or better results in all the experiments.

We are currently working on extending our approach to problems involving linear
observation operators ({\it e.g.}, blurs), other non-additive and
non-Gaussian noise models, such as Poissonian observations \cite{Figueiredo_Bioucas_SSP09},
\cite{SetzerSteidlTeuber}, and other regularizers, such as those based on frame representations.

\section*{Acknowledgments} We would like to thank the authors
of \cite{Huang09}, Yu-Mei Huang, Michael K. Ng, and You-Wei Wen, for allowing
us to use their code and original images, and to Mila Nikolova and Jalal
Fadili for providing us the images and results from \cite{DurandFadiliN08}.

\end{document}